\documentclass[article]{amsart}
\usepackage{microtype}
\usepackage{amsmath,amssymb, mathdots}
\usepackage{caption}
\usepackage{hyperref}
\usepackage[utf8]{inputenc}
\usepackage{mathtools}
\usepackage{color,graphicx}
\usepackage{thmtools, thm-restate}
\usepackage{enumitem}
\usepackage{soul}
\usepackage{stmaryrd}
\usepackage{tensor}
\reversemarginpar
\usepackage{algorithm}
\usepackage{algpseudocode}

\usepackage[colorinlistoftodos,textsize=tiny,textwidth=31mm]{todonotes}
\usepackage{etoolbox}
\usepackage{mathabx}
\makeatletter

\patchcmd{\@addmarginpar}{\ifodd\c@page}{\ifodd\c@page\@tempcnta\m@ne}{}{}
\providecommand*{\cupdot}{%
  \mathbin{%
    \mathpalette\@cupdot{}%
  }%
}
\newcommand*{\@cupdot}[2]{%
  \ooalign{%
    $\m@th#1\cup$\cr
    \hidewidth$\m@th#1\cdot$\hidewidth
  }%
}
\providecommand*{\bigcupdot}{%
  \mathop{%
    \vphantom{\bigcup}%
    \mathpalette\@bigcupdot{}%
  }%
}
\newcommand*{\@bigcupdot}[2]{%
  \ooalign{%
    $\m@th#1\bigcup$\cr
    \sbox0{$#1\bigcup$}%
    \dimen@=\ht0 %
    \advance\dimen@ by -\dp0 %
    \sbox0{\scalebox{2}{$\m@th#1\cdot$}}%
    \advance\dimen@ by -\ht0 %
    \dimen@=.5\dimen@
    \hidewidth\raise\dimen@\box0\hidewidth
  }%
}

\newcommand{\RN}[1]{%
  \textup{\uppercase\expandafter{\romannumeral#1}}%
}

\newenvironment{bcolarray}{\left[\begin{array}{@{}c@{}}}{\end{array}\right]}

\makeatother
\reversemarginpar
\newcommand{\cE}{\mathcal{E}}

\newcommand{\cO}{\mathcal{O}}

\newcommand{\cU}{\mathcal{U}}

\newcommand{\bitm}{\begin{itemize}}
\newcommand{\eitm}{\end{itemize}}
\newcommand{\bitme}{\begin{enumerate}[label=(\roman*),leftmargin=0.25in]}
\newcommand{\eitme}{\end{enumerate}}
\newcommand{\beq}{\begin{equation}}
\newcommand{\eeq}{\end{equation}}
\newcommand{\btcb}{\begin{tcolorbox}}
\newcommand{\etcb}{\end{tcolorbox}}
\def\bals#1\eals{\begin{align*} #1 \end{align*}}
\def\bal#1\eal{\begin{align} #1 \end{align}}

\newcommand{\where}{\quad \text{ where } }

\newcommand\Dom\Omega

\newcommand\RR{\mathbb{R}}

\newcommand\NN{\mathbb{N}}
\newcommand\NNz{\mathbb{N}_0}

\newcommand\ZZ{\mathbb{Z}}

\newcommand\Lap\Delta

\newcommand\abs[1]{\left\lvert #1 \right\rvert}

\def\bpde#1\epde{\[\left\{\begin{aligned}#1\end{aligned}\right. \]}
\def\inbpde#1\inepde{\left\{\begin{aligned}#1\end{aligned}\right.}
\def\binpde#1\einpde{\left\{\begin{aligned}#1\end{aligned}\right.}

\def\cU{\mathcal{U}}

\def\half{\frac{1}{2}}

\def\cU{\mathcal{U}}

\def\b0{\mathbf{0}}

\def\bbmat{\begin{bmatrix}[r]}
\def\ebmat{\end{bmatrix}}

\newcommand{\barr}{\begin{array}}
\newcommand{\ea}{\end{array}}
\newcommand{\bea}{\begin{eqnarray}}
\newcommand{\eea}{\end{eqnarray}}
\newcommand{\bt}{\begin{table}}
\newcommand{\et}{\end{table}}

\DeclareMathOperator\Id{Id}
\DeclareMathOperator\Span{span}

\DeclareMathOperator\rg{range}

\DeclareMathOperator\diag{diag}
\DeclareMathOperator\blockdiag{diag}

\newtheorem{theorem}{Theorem}[section]
\newtheorem{lemma}[theorem]{Lemma}
\newtheorem{corollary}[theorem]{Corollary}
\newtheorem{definition}[theorem]{Definition}

\numberwithin{equation}{section}

\newcommand\tturl[1]{{\tt \scriptsize [\url{{#1}}]}}

\DeclarePairedDelimiter\dparen{\lparen}{\rparen}
\DeclarePairedDelimiter\dbrace{\lbrace}{\rbrace}

\newcommand{\remove}[1]{\textcolor{red}{[Text removed.]}}%

\graphicspath{{../figures/}{../python/}}

\newcommand\brem{\mathbin{\%}}      % remainder, as binary operation

\newcommand\iu{\mathrm{i}}      % imaginary unit
\newcommand\nm[1]{{#1}^{n}_{m}}

\newcommand\nmz[2]{{#1}^{n}_{#2}}
\newcommand\nmzT[2]{{#1}^{n \top}_{#2}}
\newcommand\nmzD[2]{{#1}^{n \dagger}_{#2}}

\newcommand\nR{R^n}

\newcommand\npmzU[1]{\nmz{U}{(#1)}}

\newcommand\npmzVT[1]{\nmzT{V}{(#1)}}
\newcommand\npmzV[1]{\nmz{V}{(#1)}}

\newcommand\npmzSig[1]{\nmz{\Sigma}{(#1)}}

\newcommand\npmzS[1]{\nmz{S}{(#1)}}

\newcommand\npmzSD[1]{\nmzD{S}{(#1)}}

\begin{document}

\ifpdf
  \DeclareGraphicsExtensions{.pdf, .jpg, .tif}
\else
  \DeclareGraphicsExtensions{.eps, .jpg}
\fi

\title{An explicit spectral decomposition of the ADRT}

\author{Weilin Li}
\address{Department of Mathematics, 
  City University of New York, City College and Graduate Center,
  New York, NY 10031}
\email{{\tt wli6@ccny.cuny.edu}}

\author{Karl Otness}
\address{Courant Institute of Mathematical Sciences, %
    New York University, New York, NY 10012}%
\email{{\tt karl.otness@nyu.edu}}%

\author{Kui Ren}
\address{Department of Applied Physics and Applied Mathematics, Columbia
    University, New York, NY 10027}
\email{{\tt kr2002@columbia.edu}}

\author{Donsub Rim}
\address{Department of Mathematics, %
    Washington University in St. Louis, St. Louis, MO 63130} %
\email{{\tt rim@wustl.edu}}%

\maketitle

\begin{abstract}
  The approximate discrete Radon transform (ADRT) is a hierarchical multiscale
  approximation of the Radon transform. In this paper, we factor the ADRT into a
  product of linear transforms that resemble convolutions and derive an explicit
  spectral decomposition of each factor. We further show that this implies---for
  data lying in the range of the ADRT---that the transform of an $N \times N$
  image can be formally inverted with complexity $\cO(N^2 \log^2 N)$. We
  numerically test the accuracy of the inverse on images of moderate size and
  find that it is competitive with existing iterative algorithms in this special
  regime.
\end{abstract}

\section{Introduction} \label{sec:intro}

The continuous Radon transform takes a function defined over the Euclidean
plane and produces a new function over the space of all straight lines,
assigning per line the value of the corresponding line integral of the original
function \cite{Helgason2010,Natterer}. 
The approximate discrete Radon transform (ADRT)
\cite{Brady98adrt,Gotz96fdrt,Press06drt} is a fast algorithm for computing a
multiscale discretization of the continuous Radon transform.
ADRT is applicable in various contexts, such as in computerized tomography and
inverse problems \cite{Natterer,Epstein2007}, solutions to partial differential
equations \cite{Folland2020,rim18split}, and in sliced Wasserstein barycenters
or interpolants \cite{Bonneel15slice,Rim18mr}, to name a few.

Given an $N \times
N$ image with $N = 2^n$ for some positive integer $n$, the ADRT $\nR$ transforms
that image in $\RR^{N \times N}$ to one in $\RR^{4 \times N \times (\frac{3}{2}
N - \frac{1}{2})}$ using $\cO(N^2 \log N)$ floating point operations 
 where $\log$ refers to the logarithm of base two throughout,
and this
transformed data is an approximation of the continuous Radon transform of that
image. 
The ADRT is defined as a collection of
sums over the discrete values $x(i,j)$ of an image $x \in \RR^{N \times N}$.
The digital lines $D_n$ are a subset
of the indices $(i,j)$ that will be summed together; that is, $R^n$ produces a
real value for each digital line $D_n$. The digital line is parametrized by two
integral indices $(h, s)$ belonging to a subset of $\ZZ \times \ZZ$. 
The \emph{single-quadrant ADRT} $R^n_n : \RR^{N \times N} \to \RR^{N \times
(\frac{3}{2}N - \frac{1}{2})}$ of an image $x \in \RR^{N \times N}$ is given by
\begin{equation} 
    R^n_n [x](h, s) := \sum_{(i,j) \in D_n(h,s)} x(i,j)
\end{equation}
The digital lines $D_n(h, s)$ have a recursive structure: $D_n(h,s)$ of length
$N = 2^n$ is a union of two digital lines of half-lengths $N/2 = 2^{n-1}$
by the relations
\begin{equation}
  \begin{aligned}
  D_n(h, 2s) &= D_{n-1}^{(L)} (h, s) \cup D_{n-1}^{(R)}(h+s, s),
  \\
  D_n(h, 2s+1) &= D_{n-1}^{(L)} (h, s) \cup D_{n-1}^{(R)}(h+s+1, s).
  \end{aligned}
\end{equation}
This hierarchical structure allows one to compute $R^n_n$ in $\cO(N^2 \log
N)$ operations. Since the sum $R^n_n$ can now be computed from two ADRTs
$R^n_{n-1}$ of left and right half-images $x^{(L)}$ and
$x^{(R)}$ both in $\RR^{N \times N/2}$.  Correspondingly, at each level $m = 1,
... , n$ two $R^n_{m-1}$ are summed to form $R^n_{m}$,
\begin{equation}\label{eq:Rnm_split}
  \begin{aligned}
   R^n_m [x]  (h, 2s) &= R^n_{m-1} [x^{(L)}] (h, s)
                            +
                            R^{n}_{m-1} [x^{(R)}] (h+s, s),
                            \\
   R^n_m  [x] (h, 2s+1) &= R^{n}_{m-1} [x^{(L)}] (h, s)
                            +
                           R^{n}_{m-1} [x^{(R)}] (h+s+1, s).
  \end{aligned}
\end{equation}
 Thus, after such sums over $m = 1, ...\,, n$ each costing
$\cO(N^2)$ operations, $R^n_n$ is produced.

\begin{figure}
  \centering
  \begin{tabular}{ll}
    {\scriptsize (a)} & {\scriptsize (b)} \\
  \begin{minipage}{0.3\textwidth}
  \includegraphics[width=0.9\textwidth]{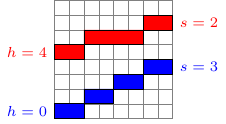}
  \\ ~~ \\ ~~ \\ ~~
  \end{minipage}
  &
  \begin{minipage}{0.5\textwidth}
  \includegraphics[width=0.85\textwidth]{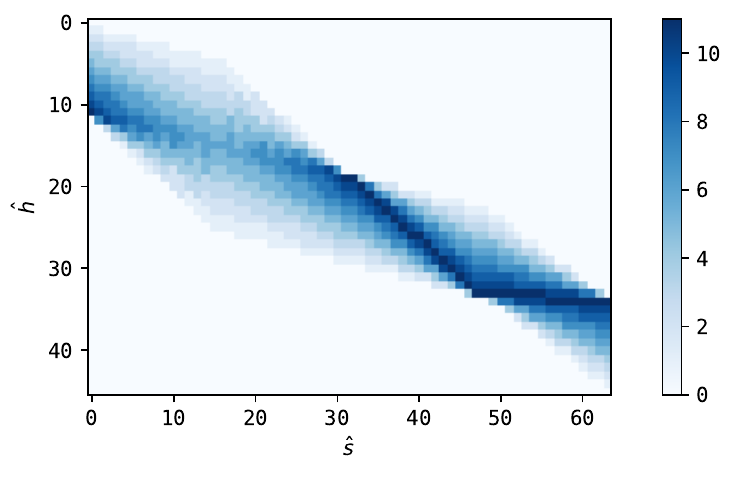}
  \end{minipage}
\\
{\scriptsize (c)}
\\
\multicolumn{2}{l}{
  \begin{minipage}{0.75\textwidth}
  \includegraphics[width=0.9\textwidth]{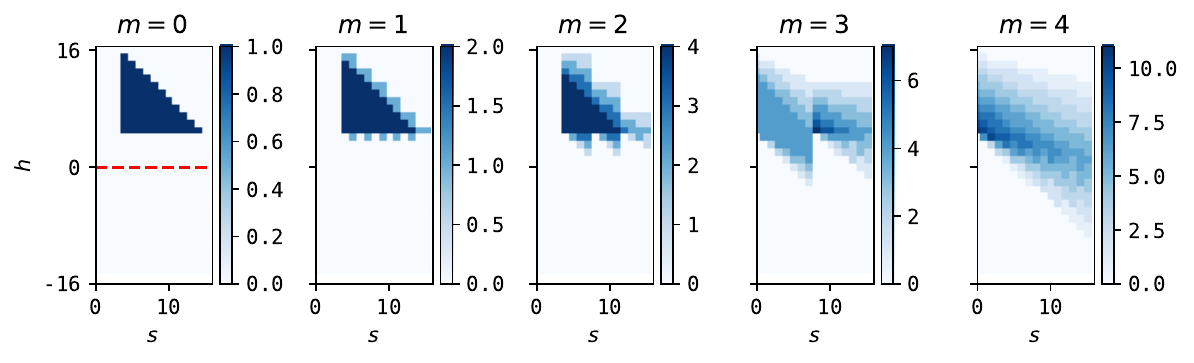}
  \end{minipage}
}
  \end{tabular}
  \caption{  (a) A depiction of two digital lines in $8 \times 8$ images.
  $D_8(4,2)$ and $D_8(0,3)$ are shown,
  (b) full ADRT plot stitching four single-quadrant ADRTs applied to
  rotated and flipped images,
  and
  (c) An example of the single-quadrant
  ADRT of an image taking on values of $0$s and $1$s. The red dashed line
  denotes the extent of the square image $u \in \RR^{16 \times 16}$. $R^n_{(m)}
  [u]$ is shown for $m = 0, 1, ...\,, 4$. One single-quadrant
  ADRT is illustrated, for lines with the angles in $[-\frac{\pi}{2},
  -\frac{\pi}{4}]$.}
  \label{fig:adrt_part}
\end{figure}

The digital line $D_n$ is
naturally defined for lines corresponding to angles in $[0, \frac{\pi}{4}]$, and
so the above computation for these angles is called the single-quadrant
ADRT \cite{Rim20iadrt}.
Fig.~\ref{fig:adrt_part} (a) shows a depiction of two digital lines.
A plot showing each level of a single-quadrant ADRT
is shown in Fig.~\ref{fig:adrt_part} (b): each pixel in the subsequent level is a
sum of two pixels from the previous level.

The \emph{full ADRT} (or simply \emph{the ADRT}) is computed by rotating and
flipping the original image to obtain the approximation to line integrals
for all angles in $[-\frac{\pi}{2}, \frac{\pi}{2}]$. The four single-quadrant ADRTs
of the image can be stitched together in a contiguous image, as shown in
Fig.~\ref{fig:adrt_part} (c) where a very low resolution was chosen
for illustration.
The single-quadrant ADRT and the full ADRT will be denoted by
$R^n_{n}$ and $R^n = R^n_{(n)}$, respectively. 

In this work, we derive a new spectral decomposition of the ADRT, and
demonstrate that the decomposition leads to a fast and explicit computation of
the inverse of the ADRT in the setting where the data lies in its range.

As we will illustrate below, the ADRT can be factored as
\begin{equation} \label{eq:Rn_factor}
  \nR = \npmzS{n} \npmzS{n-1} \cdots \npmzS{1}
\end{equation}
where $\npmzS{m}$ for $m = 1, ...\,,n,$ represent summation operations \eqref{eq:Rnm_split} at each
multiscale level of the ADRT corresponding to the length scale $\sim 2^{-(n-m)}$.
The main contribution of this work is the derivation of an explicit spectral
decomposition of each of the factors,
\begin{equation} \label{eq:Sm_decomp}
  \npmzS{m} = \npmzU{m} \npmzSig{m} \npmzVT{m},
  \qquad
  m = 1, 2, ...\,,n,
\end{equation} 
which possesses certain desirable properties: (1) The matrices $\npmzU{m}$ and
$\npmzV{m}$ are either orthogonal or contain orthogonal square subblocks of
maximal size; (2) $\npmzSig{m}$ is diagonal with positive entries; and (3) the
columns of $\npmzU{m}$ and $\npmzV{m}$ are each made up of columns of discrete
sine transforms (DSTs) or discrete cosine transforms (DCTs), and the
multiplication by these matrices or their transposes can be performed in 
$\cO(N^2 \log N)$ operations. 

We use a notation that resembles that for the singular value decomposition (SVD)
for this decomposition since it is indeed the SVD of $S^n_{m}$ for the levels $m
> 1$, and it is a spectral decomposition with an SVD-like structure for the
level $m = 1$. We will remind the reader of this distinction throughout.
  
The fact that discrete sine and cosine transforms make up these columns is not a
complete surprise, considering how the Fourier slice theorem \cite{Natterer}
relates the continuous Radon transform to the Fourier transform in a
straightforward manner.  However, in the discretized setting the image and the
transformed data are represented as pixelated images on Cartesian grids, and as
a result the precise relation between a discretized Radon transform and the
discretized Fourier transform is not straightforward to pin down. The relation
depends delicately on the choice of discretization; for example, the
preservation of such a relation is an important consideration in the derivation
of the digital Radon transform via pseudo-polar Fourier transforms
\cite{Averbuch08ppft,Averbuch08drt,Averbuch20033d}. The explicit
decomposition we present here \eqref{eq:Sm_decomp} establishes a precise discrete
relation for the ADRT: since the columns are made up of DSTs or DCTs, one can
equate the discrete transforms on the ADRT data to a variant of the discrete
Fourier transform in the image domain. That is, our spectral decomposition
implies a fully discrete version of the Fourier slice theorem that is specific
to the ADRT.
  
The spectral decomposition is useful for inverting the ADRT.  Given some ADRT
data $b \in \RR^{4 \times N \times (\frac{3}{2} N - \frac{1}{2})}$ lying in the
range $b \in \rg(R^n)$, we consider the problem of finding $x \in \RR^{N
\times N}$ satisfying
\begin{equation} \label{eq:prob}
  \nR [x] = b.
\end{equation}
Broadly speaking, most readers will be more interested in the general
inverse problem $\min_x \lVert  \nR[x] - b \rVert$ for some norm $\lVert \cdot
\rVert$, without the assumption that $b \in \rg(R^n)$.  While that general
problem
is outside the scope of this paper, it is important to note that the range
characterization result \cite{LRR23} showed $\rg(R^n)$ to be highly
structured. Exploiting this structure, it is potentially possible to numerically
project the
data $b$ directly to $\rg(R^n)$, and the projection will reduce the general
problem to our setting. Hence, we believe the results in this
paper will be relevant to the general problem as well. 

The decomposition \eqref{eq:Sm_decomp} yields an explicit expression for the
pseudo-inverse for each $\npmzS{m}$.  With the pseudo-inverses $\npmzSD{m}$ of
individual levels $m = 1, ...\,,n$ at hand, the inverse $x$ satisfying
\eqref{eq:prob} is given by
\begin{equation} \label{eq:pinv_sol}
  x
  =
  \npmzSD{1}
  \cdots
  \npmzSD{n-1}
  \npmzSD{n} b.
\end{equation}
Due to the block structures of $\npmzU{m}$ and $\npmzV{m}$ and due to their
columns being expressed in terms of DSTs and DCTs, this pseudo-inverse can be
computed in $\cO(N^2 \log N)$ serial operations, or in $\cO(N \lceil \frac{N}{P}
\rceil \log N)$ operations with $P < N$ parallel processors.  We will call this
fast and explicit computation of the inverse, the \emph{Spectral Pseudo-Inverse,
Fast and Explicit} (SPIFE) of the ADRT.
  
SPIFE is exact if $b \in \rg(\nR)$ and if there were no errors in the floating
point arithmetic. With finite-precision floating point arithmetic, the accuracy
of SPIFE depends on the precision and the size of the image.  For
double-precision and image sizes that are moderately small $N \le 2^8$, this
explicit inverse outperforms existing iterative inverses, the generic conjugate
gradient (CG) algorithm \cite{Gre97} and the full multi-grid (FMG) inverse
devised for the ADRT \cite{Press06drt}, in terms of accuracy per
workload for the case. We will present numerical experiments demonstrating this
performance.
  
We perform additional numerical experiments and explain why finite-precision
SPIFE loses accuracy. As the pseudo-inverses are applied from right to left in
\eqref{eq:pinv_sol}, the floating point calculation errors accumulate in the
levels $m > 1$ due to the large values in $(\Sigma^n_m)^{-1}$ appearing inside
the pseudo-inverse $\npmzSD{m}$ and significant digits are progressively lost.
An accurate image is recovered only if sufficient digits are retained before
reaching level $m=1$, at which point $\npmzSD{1}$ stabilizes the inverse by removing some
of the unstable components.  For a given precision level, this means the image sizes must not be too large. 

A closely related fact is that the explicit inverse of cost $\cO(N^2 \log N)$
for the single-quadrant ADRT is highly unstable
\cite{Rim20iadrt,Press06drt,OliviaGarcia2022}.  This explicit single-quadrant
inverse is---asymptotically as the resolution $N$
increases---related to the limited-angle problem \cite{LouisRieder89}%
\cite[Section~VI.2]{Natterer}
for the continuous Radon transform: that
is, for a sequence of images $x_k \in \RR^{N_k \times N_k}$ and data $b_k \in
\RR^{4 \times N_k \times (\frac{3}{2} N_k - \frac{1}{2})}$ that approaches a
function and its Radon transform as $k \to \infty$, respectively, the inversion
of the single-quadrant ADRT approaches a limited-angle problem.   Due to this
connection, the instability is inherent, and one does not generally expect to
find numerical techniques for controlling it. In contrast, SPIFE is not
inherently unstable, as it utilizes ADRT data from all four quadrants, and so
asymptotically its stability should match that of the (full-angle) inversion of
the continuous Radon transform; this inversion is only mildly ill-conditioned
\cite{Natterer}. In this sense, one may view SPIFE as a stable extension of the
na\"ive explicit single-quadrant inverse to the full ``all-quadrants'' ADRT. We
provide a comparison in our experiments.

Previously, Press observed in his experimentation with the FMG inverse a certain
scaling between the number of iterations and error reduction, and claimed that
his FMG inverse would achieve a desired error threshold with $\cO(N^2 \log^3 N)$
work  \cite{Press06drt}. Based on this observed scaling, he further conjectured
that an iterative algorithm of $\cO(N^2 \log^q N)$ with $2 \le q \le 3$ could
exist and that $q=2$ would be optimal.  These complexity bounds were never
rigorously proved, however. As for SPIFE, it has precisely the best
conjectured complexity $\cO(N^2 \log^2 N)$, and if (1)
the provided data $b$ lies in the
range up to double precision, and (2) the image size is not large ($N \le 2^8$),
then the explicit pseudo-inverse has competitive accuracy.
Hence, our
result marks partial progress towards the conjecture of Press.
  
However, as mentioned above, the key practical motivation for studying the
ADRT lies in the search
for a faster inverse postulated in the unrestricted conjecture: an optimally fast
Moore-Penrose inverse that is accurate up to the inherent conditioning of $R^n$
regardless of the image size or whether the data lies in the range.  The results
here offer promise towards that goal. To achieve it, SPIFE must be aided by two
further developments: (1) a fast projection algorithm that is able to
orthogonally project the data to the range of the ADRT \cite{LRR23} to numerical
precision,
and (2) a stabilization technique that prevents numerical errors from
accumulating when applying the pseudo-inverse. These topics will be pursued in
future works.

\section{Notations, preliminaries, and outline} 

We recall a few definitions and notations from \cite{Rim20iadrt,LRR23}
and then outline the key results of the paper.

\subsection{Notation and preliminaries}

Let $\ZZ$ denote the integers, $\NNz$ the non-negative integers, and $\NN$ the
natural numbers.  Define a set of indices $[m] := \{i \in \NNz \mid i < m\}$ for
$m \in \NNz$, and $[m]_1 := \{i + 1 \mid i \in [m]\}$. For simplicity of
notation, we will use members of $[m] \times [n]$ for row and column indexing of
a matrix in $\RR^{m \times n}$ in lieu of $[m]_1 \times [n]_1$ so that the
indexing will start from $0$ rather than $1$.  $\Id^n$ will denote the identity
in $\RR^{2^n \times 2^n}$ whereas $\Id_m$ denotes the identity in $\RR^{m \times m}$.
Let us define the remainder $\brem: \NNz \times \NNz \to \NNz$ given by $a \brem
b := a - \lfloor a/b \rfloor \cdot b$ if $b \neq 0$ and by $a \brem b := a$ if
$b = 0$.

To describe ADRT data, it is convenient to allow tuples of indices with negative
values for each index, and to allow variable bounds for the index range. 

We define the index set $E^n_m$ as
\begin{equation}
  E^n_m
  =
  \dbrace[\big]{
  (h, s) \in \ZZ \times [2^n]
  \mid
    - (s \brem 2^m) \le h < 2^n
  }.
\end{equation}
We partition $E^n_{m}$ into subsets $E^n_{m, \ell}$, 
\begin{equation}
  E^n_m
  =
  \bigcup_{\ell \in [2^{n-m}]} E^n_{m, \ell},
\end{equation}
where $E^n_{m, \ell}$ is the restriction of the second index to the range
$\ell 2^m + [2^m]$, 
\begin{equation}
  E^n_{m, \ell}
  =
  \dbrace[\big]{
  (h, s) \in \ZZ \times [2^n]
  \mid
    -(s \brem 2^m) \le h < 2^n,
    s \in \ell 2^m + [2^m]
  }.
\end{equation}
We further partition $E^n_{m, \ell}$ into subsets $E^n_{m, \ell, t}$
\begin{equation}
  E^n_{m, \ell}
  =
  \bigcup_{t \in [2^{m}]} E^n_{m, \ell, t},
\end{equation}
where $E^n_{m, \ell, t}$ is the restriction of the second index to a single
value 
\begin{equation}
  E^n_{m, \ell, t}
  =
  \left\{
  (h, s) \in \ZZ \times [2^n]
  \mid
    -(t \brem 2^m) \le h < 2^n,
    s = \ell 2^m + t
  \right\}.
\end{equation}

\begin{figure}
\centering
\includegraphics[width=0.45\textwidth]{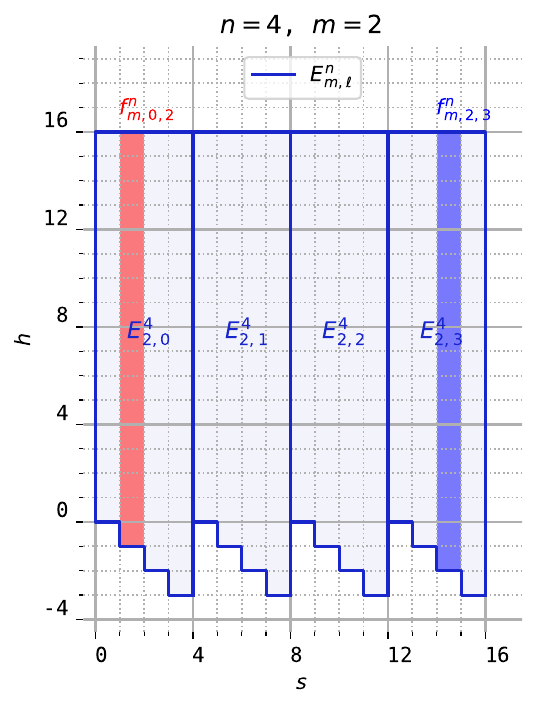}
\caption{An example of the vectors $f^n_{m, \ell, t}$, illustrating two different instances of these vectors and labeling the index sets for each of the four sections of the ADRT.}
\label{fig:fnmlt}
\end{figure}

Next, we define functions $f \in \RR^{E^n_n}$ on these index sets. Such
a function can also be viewed as a matrix with suitable zero padding, or an
array of dimension two with variable column sizes. Let us define the restriction
$f^n_{m, \ell} \in \RR^{E^n_{m, 0}}$ of $f \in \RR^{E^n_n}$ as 
\begin{equation}
    f^n_{m, \ell}(h, s) := f(h, s + \ell 2^m),
    \quad (h, s) \in \RR^{E^n_{m, 0}},
\end{equation}
and similarly the restriction $f^n_{m, \ell, t} \in \RR^{E^n_{m, 0, t}}$ of
$f^n_{m, \ell} \in \RR^{E^n_{m, 0}}$ as 
\begin{equation}
    f^n_{m, \ell, t}(h)
    := f^n_{m, \ell} (h, t)
    = f(h, t + \ell 2^m),
    \quad (h, t) \in \RR^{E^n_{m, 0, t}}.
\end{equation}
See Fig.~\ref{fig:fnmlt} for an illustrative example.

The single-quadrant ADRT at the individual level $m$ is denoted by
\begin{equation} \label{eq:sq_Snm}
  \nm{S}: \RR^{E^n_{m-1}} \to \RR^{E^n_m},
  \qquad
  S^n_m[f] = g,
\end{equation}
where for each $(h, s) \in E^n_{m-1, 0}$, we let
\begin{equation}
  \begin{aligned}
    g^n_{m, \ell} (h, 2s + 1)
    &:=
    f^n_{m, 2\ell} (h, s) + f^n_{m, 2\ell + 1} (h+s+1, s),
    \\
    g^n_{m, \ell} (h, 2s)
    &:=
    f^n_{m, 2\ell} (h, s) + f^n_{2\ell+1}(h+s, s).
  \end{aligned}
\end{equation}

Now we set additional notations and definitions related to the full ADRT, as
opposed to those related to the single-quadrant ADRT considered thus far. We
distinguish the analogous mappings by putting a parenthesis around the level
index $m$. For example, $\nm{S}$ denotes a single-quadrant ADRT level whereas
$S^n_{(m)}$ denotes the full ADRT level.

\begin{definition} \label{def:Snm}
The level $m$ of the full ADRT for image $f \in \RR^{2^n \times 2^n}$ is denoted
by $S^n_{(m)}$, which we define separately for the case $m = 1$ and the case $m =
2, 3, ...\,, n$.
\begin{itemize}

  \item For $m = 1$, we define $S^n_{(m)}: \RR^{2^n \times 2^n} \to \RR^{4
  \times E^n_1}$ and let
  \begin{equation} \label{eq:Sn1}
    S^n_{(1)}
    :=
    \left[
      S^n_1 T^n_{\RN{1}} ,
      S^n_1 T^n_{\RN{2}} ,
      S^n_1 T^n_{\RN{3}} ,
      S^n_1 T^n_{\RN{4}} 
    \right]^\top,
  \end{equation}
  where $T^n_{\RN{1}}, ...\,, T^n_{\RN{4}}: \RR^{2^n \times 2^n} \to \RR^{2^n
  \times 2^n}$ are permutations corresponding to the four quadrants $\RN{1},
  ...\,, \RN{4}$ given by
  \begin{equation} \label{eq:permute}
    \left\{
    \begin{aligned}
      (T^n_{\RN{1}} f)(i, j) & = f(j \,, 2^n-1-i),
      \\
      (T^n_{\RN{2}} f)(i, j) & = f(2^n-1-i \,, j),
      \\
      (T^n_{\RN{3}} f)(i, j) & = f(i \,, j),
      \\
      (T^n_{\RN{4}} f)(i, j) & = f(2^n-1-j \,, 2^n-1-i),
    \end{aligned}
    \right.
    \quad
    i, j \in [2^n].
  \end{equation}
  We refer to the first level as the \emph{cross-quadrant level}.

  \item  For $m = 2, ... \,,n$, let $S^n_{(m)}: \RR^{4 \times E^n_{m-1}} \to
  \RR^{4 \times E^n_m}$ with 
  \begin{equation} \label{eq:Snm}
    S^n_{(m)}
    :=
    \Id_4 \otimes \nm{S},
  \end{equation}
  that is, the block diagonal matrix,
  \begin{equation} \label{eq:Snm_blockdiag}
    S^n_{(m)}
    =
    \blockdiag \left[ \nm{S} , \nm{S} ,  \nm{S},  \nm{S} \right].
  \end{equation}
  We refer to these subsequent levels as \emph{single-quadrant levels}.

  \end{itemize}
\end{definition} 

\begin{figure}
\centering
\includegraphics[width=0.8\textwidth]{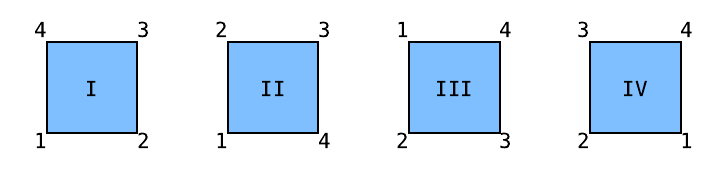}
\caption{The orientation of the four quadrants. The four corners are numbered
          to indicate the precise orientation, and the orientation of quadrant
          $\RN{3}$ matches that of the original image, that is, $T^n_\RN{3} = \Id^n$ \eqref{eq:permute}.}
\label{fig:quadrant_orient}
\end{figure}

See Fig.~\ref{fig:quadrant_orient} for a depiction of the orientation of the four
quadrants.

Finally, $R^n_{(n)}$ will denote the full ADRT, as given by $R^n_{(n)} =
S^n_{(n)} S^n_{(n-1)} \, ...  \, S^n_{(1)}$ \eqref{eq:Rn_factor}. We
sometimes omit the subscript of $R^n_{(n)}$ and write $R^n$ instead.

\subsection{Outline}

We provide an outline of this paper. As discussed in the introduction
(Sec.~\ref{sec:intro}), the main goal of this paper is to show that the individual
level $\npmzS{m}$ of the full ADRT (Def.~\ref{def:Snm}) has a spectral decomposition
in an explicit form \eqref{eq:Sm_decomp}.  We first derive the decomposition of
$\npmzS{m}$ for levels $m > 1$ \eqref{eq:Snm} in Sec.~\ref{sec:svd}, which is given
straightforwardly by the convolution theorem upon certain permutations.  The
decomposition of $\npmzS{m}$ for level $m = 1$ \eqref{eq:Sn1} is more involved
due to the presence of the rotations and reflections $T^n_{\RN{1}}, ...\,,
T^n_{\RN{4}}$ in \eqref{eq:permute}, but we show it is still possible to derive a
spectral decomposition for this case in Sec.~\ref{sec:spife} by exploiting discrete
directional patterns we devise.

An important consequence of this decomposition is that the inverse problem
\eqref{eq:prob} can be solved efficiently using the pseudo-inverses of
the $(\npmzS{m})_{m=1}^n$,  assuming the supplied data $b$ is in the range of $\nR$.  We
referred to this fast computational inverse as SPIFE in the Introduction, and it
is defined in Sec.~\ref{sec:spife_def}.  Results from numerical experiments
regarding SPIFE are reported in Sec.~\ref{sec:numerics}. Therein we provide
comparisons with existing iterative inverses and detail numerical issues
that arise when double precision floating point arithmetic is used in SPIFE.

\section{SVD of single-quadrant levels} \label{sec:svd}

In this section, we derive the explicit SVDs of each level of the ADRT $S^n_{(m)}$
for levels $m > 1$. The result follows from the simple observation that, upon
certain permutations, $S^n_{(m)}$ can be written as a block diagonal matrix, and
the individual block is either a matrix representing a convolution or a
transposed convolution. For convolution and transposed
convolutions, one can derive an explicit SVD where the singular vectors are
discrete Fourier bases. 

Let us define for a given $t \in \NN$ the convolution matrices $K_t^\pm \in
\RR^{(t+1) \times t}$ 
\begin{equation} \label{eq:Kt}
  K_t^\pm
  :=
  \begin{bmatrix}
  \pm 1 &   &        &        &   \\
    1 &\pm 1 &        &        &   \\
      & 1 &\pm 1      &        &   \\
      &   & \ddots & \ddots &   \\
      &   &        & 1      &\pm 1 \\
      &   &        &        & 1 \\
  \end{bmatrix}.
\end{equation}
In common machine learning terminology, this is a discrete convolution
operator with kernels $[1, \pm 1]^\top$, stride one, and zero
padding one; for a nice visual guide to these terms, see
\cite{dumoulin2018}. 
The SVD of $K_t^\pm$ is given explicitly as a consequence of the convolution
theorem (see standard texts \cite[Sec.~5.3.]{Proakis1996} or \cite[Thm.~3.7.]{Mallat2009}). Details are given in the following
lemma.

\begin{lemma} \label{lem:conv}
The SVD of $K_t^\pm$  (\/$t \in \NN$\/) is given by
\begin{equation}
  K_t^\pm = U_t^\pm \Sigma_t^\pm V_t^{\pm \top}
\end{equation}
where $\Sigma_t^\pm$ is a diagonal matrix whose entries are given by
\begin{equation} \label{eq:Ktsk}
  \sigma_k^+
  =
  2 \cos \left( \frac{k \pi}{2 t + 2} \right),
  \quad
  \sigma_k^-
  =
  2 \sin \left( \frac{k \pi}{2 t + 2} \right),
  \quad
  k \in [t]_1,
\end{equation}
the columns of $U_t^\pm$ are
\begin{equation} \label{eq:Ktuk}
    (u_{k}^+)_j
    =
    \sin \left( \frac{k \pi (j + \half)}{t + 1}  \right),
    \quad
    (u_{k}^-)_j
    =
    \cos \left( \frac{k \pi (j + \half)}{t + 1}  \right),
    \quad
    j \in [t],
\end{equation}
and the columns of $V_t^\pm$ are
\begin{equation} \label{eq:Ktvk}
  \begin{aligned}
    (v_k^\pm)_j
    &=
    \sin \left( \frac{k \pi (j+1)  }{t+1} \right),
  \end{aligned}
  \quad
  j \in [t].
\end{equation}
\end{lemma}

\begin{proof}

We only prove the statements for $K^+_t$; the case for $K^-_t$
is very similar. Let us assume the eigenvectors of $K_t^+ K_t^{+ \top}$ are of the
form $(w_q)_j = e^{\iu q(j+1)}$. Then $(w_q)_{j=0}^{t}$ satisfying $K_t^+ K_t^{+
\top} w_q = \lambda w_q$ is equivalent to $(w_q)_{j=-1}^{t+1}$ satisfying
\begin{equation} \label{eq:midrows}
  (w_q)_{j-1} + 2 (w_q)_j + (w_q)_{j+1} = \lambda (w_q)_j,
\end{equation}
together with
\begin{equation} \label{eq:firstlastrows}
    (w_q)_{-1} + (w_q)_{0} = 0,
    \quad
    (w_q)_{t} + (w_q)_{t+1} = 0.
\end{equation}

We see that \eqref{eq:midrows} is satisfied for $w_q$ with eigenvalues 
$\lambda = 4 \cos^2(q/2)$ since
\begin{equation} 
  \begin{aligned}
    (w_q)_{j-1} +  2 (w_q)_j + (w_q)_{j+1}                          
                  & = (e^{\iu q(j-1)} + 2 e^{\iu qj} + e^{\iu q(j+1)}) 
                  e^{\iu q} \\
                  & = (e^{-\iu q} + 2 + e^{\iu q})e^{\iu q(j+1)}         \\
                  & = 2(\cos(q) + 1) e^{\iu q (j+1)}
                    = 4 \cos^2(q/2) (w_q)_j.
  \end{aligned}
\end{equation}

As $\cos(\cdot)$ is an even function, the eigenvalues are the same for
$w_q$ and $w_{-q}$, so we can satisfy \eqref{eq:firstlastrows} by letting 
\begin{equation}
  u^+_{q} := \frac{1}{2 \iu} (w_q e^{-\iu \frac{q}{2}}
                         - w_{-q} e^{\iu \frac{q}{2}})
\end{equation}
then observing that such $v_q$ is an eigenfunction if and only if $q (t + 1) =
k \pi$ for some $k \in \ZZ$. So, the eigenfunctions are as in \eqref{eq:Ktvk}.
 The eigenvalues are
\begin{equation}
  \lambda^+_k = 4 \cos^2 \left( \frac{k \pi}{2 t + 2} \right),
  \quad
  k \in [t+1]_1,
\end{equation}
and one eigenvalue $\lambda_k$ for the case $k=t+1$ is zero. Discarding the zero
value, singular values $\sigma_k^+ = (\lambda_k^+)^{\half}$ are as in
\eqref{eq:Ktsk}.  Now,
\begin{equation}
  \begin{aligned}
  (\sigma^+_k v^+_k)_j
  &=
  \dparen[\big]{(K^+_t)^\top u^+_k}_j
  =
  \sin \left( \frac{k \pi (j + \half)}{t + 1} \right)
  +
  \sin \left( \frac{k \pi (j + \frac{3}{2})}{t + 1} \right)
  \\
  &=
  2\cos \left( \frac{k \pi}{2(t+1)} \right)
  \sin \left( \frac{k \pi (j+1)}{t + 1} \right)
  \end{aligned}
\end{equation}
for $j \in [t]$ and $k \in [t]_1$, so we have \eqref{eq:Ktuk}.

\end{proof}

\begin{figure}
  \centering
  \begin{tabular}{cc}
    \begin{minipage}{0.45\textwidth}
    \includegraphics[width=1.0\textwidth]{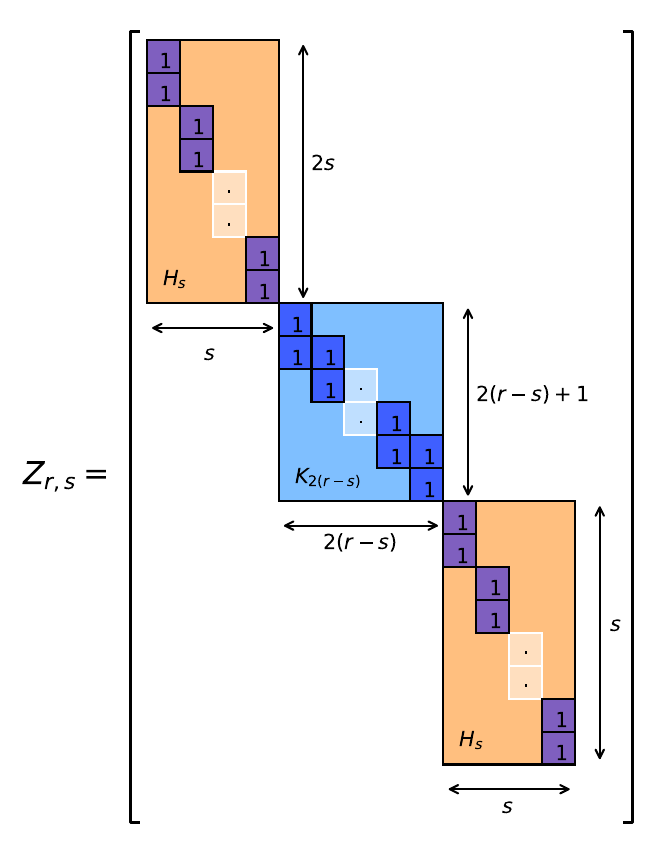}
    \end{minipage}
    &
    \begin{minipage}{0.45\textwidth}
    \includegraphics[width=1.0\textwidth]{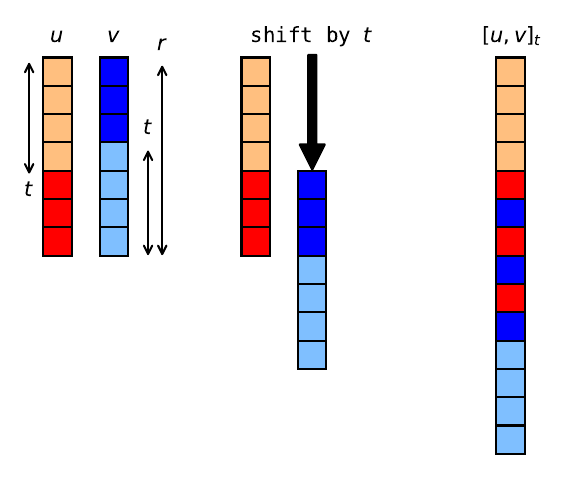}
    \end{minipage}
    \\
    (a)
    &
    (b)
  \end{tabular}
  \caption{(a) The blocks of matrix $Z_{r,s}$ \eqref{eq:Zrs}, and
  (b) an illustration of $[u, v]_t$ for $u, v \in \RR^r$ according to
  Def.~\ref{def:shiftconcat}.
  }
  \label{fig:matrix_structure}
\end{figure}

Now we make preparations to derive the block convolution form of $S^n_m$.
First, let us define $H_t \in \RR^{2t \times t}$ ($t \in \NN$) to be
\begin{equation}
  H_t
  :=
  \begin{bmatrix}
    1 &   &        &   &   \\
    1 &   &        &   &   \\
      & 1 &        &   &   \\
      & 1 &        &   &   \\
      &   & \ddots &   &   \\
      &   &        & 1 &   \\
      &   &        & 1 &   \\
      &   &        &   & 1 \\
      &   &        &   & 1 \\
  \end{bmatrix}.
\end{equation}
This is a transposed convolution matrix with a kernel of size two
and stride two \cite{dumoulin2018}, and the matrix has a trivial SVD since it has orthogonal
columns.
Then define the block matrix $Z_{r, s} \in \RR^{(2r+2s + 1) \times 2r }$
\begin{equation} \label{eq:Zrs}
  Z_{r, s}
  :=
  \diag[H_{s}, K_{2(r-s)}^{+},H_{s}].
\end{equation}
See Fig.~\ref{fig:matrix_structure} for a sketch of this matrix.

\begin{corollary}
Up to the ordering of the singular values, the SVD of $Z_{r,s}$ is 
\begin{equation}
  Z_{r,s}
  =
  U_{r,s} \Sigma_{r,s} V_{r,s}^\top
\end{equation}
where $U_{r,s}$, $\Sigma_{r,s}$, and $V_{r,s}$ are block diagonal matrices
\begin{equation}
  \begin{aligned}
  U_{r,s}
  &=
  \diag[H_s, U_{2(r-s)}, H_s],
  \\
  \Sigma_{r,s}
  &=
  \diag[\Id_s, \Sigma_{2(r-s)}, \Id_s],
  \\
  V_{r,s}
  &=
  \diag[\Id_s, V_{2(r-s)}, \Id_s].
  \end{aligned}
\end{equation}
\end{corollary}

\begin{proof}
Follows immediately due to Lem.~\ref{lem:conv} and the remarks above.
\end{proof}

Let us define
\begin{equation} \label{eq:fSnm}
  Z^n_{m, 0}
  :=
  \blockdiag
  \underbrace{ [Z_{N, 0} , Z_{N,1} , \cdots , Z_{N, 2^{m-1}-1}]}_{2^{m} \text{ blocks}}.
\end{equation}
Lem.~\ref{lem:conv} implies that
\begin{equation*}
  \begin{aligned}
    Z^n_{m, 0}
    &=
    \blockdiag [ Z_{2^n, 0}\,, Z_{2^n,1}\,, \cdots \,, Z_{2^n, 2^m-1}]
    \\
      & =
    \blockdiag[
      U_{2^n,0} \Sigma_{2^n,0} V_{2^n,0}^\top, \cdots , U_{2^n,2^m-1} \Sigma_{2^n,2^m-1} V_{2^n,2^m-1}^\top ]
      =
    U^n_{m, 0} \Sigma^n_{m, 0} V^{n \top}_{m, 0},
  \end{aligned}
\end{equation*}
where
\begin{equation}
  \begin{aligned}
    U^n_{m, 0}
      & =
      \blockdiag[ U_{2^n,0} \,,  \cdots \,, U_{2^n,2^m-1}],
    \\
    \Sigma^n_{m, 0}
      & =
      \blockdiag [\Sigma_{2^n,0} \,, \cdots \,, \Sigma_{2^n,2^m-1}],
    \\
    V^n_{m, 0}
      & =
      \blockdiag [V_{2^n,0}, \cdots, V_{2^n,2^m-1}].
  \end{aligned}
\end{equation}

As a result, if we let
\begin{equation} \label{eq:Znm}
  Z^n_m
  :=
  \Id^{n-m-1} \otimes \, Z^n_{m, 0}
  =
  \diag
  \underbrace{[
      Z^n_{m, 0}, Z^n_{m, 0}, \cdots, Z^n_{m, 0}]}_{2^{n-m-1} \text{ blocks}},
\end{equation}
and let
\begin{equation}
    U^n_m := \Id^{n-m-1} \otimes \, U^n_{m, 0},
    \quad
    \Sigma^n_m := \Id^{n-m-1} \otimes \, \Sigma^n_{m, 0},
    \quad
    V^n_m := \Id^{n-m-1} \otimes \, V^n_{m, 0},
\end{equation}
we obtain the SVD
\begin{equation} \label{eq:Znm_svd}
  Z^n_m = U^n_m \Sigma^n_m V^{n \top}_m.
\end{equation}

Next, we will describe the permutations that will transform $S^n_{(m)}$
into convolutions.  It will be handy to set a notation for a certain
alternating concatenation of two vectors with shift.

\begin{definition} [Alternating concatenation with shift]
\label{def:shiftconcat}
Given two vectors $u, v \in \RR^r$, define a vector $[ u \mid v ]_t \in \RR^{2
r}$ with $t \in [r]$ as given by 
\begin{equation}
  \begin{aligned}
    \left( [ u \mid v ]_{t} \right)_{j}
      & :=
    \begin{cases}
      u_j
      &
      \text{ if } j \in [t],
      \\
      u_{t + \lfloor (j-t-1) / 2 \rfloor}
      &
      \text{ if } j-t \in [2 r - t] \text{ and } j-t \text{ even},
      \\
      v_{\lfloor (j-t) / 2 \rfloor}
      &
      \text{ if } j-t \in [2 r - t] \text{ and } j-t \text{ odd},
      \\
      v_{j - r}
      &
      \text{ if } j - 2r + t \in [t].
    \end{cases}
  \end{aligned}
\end{equation}
The mapping $[\,\cdot \mid \cdot\,]_t :\RR^r \times \RR^r \to \RR^{2r}$ is
invertible, so we will denote the inverse operation $[\,\cdot\,]^{-1} : \RR^{2r}
\to \RR^r \times \RR^r$.
\end{definition}

As an example, given two vectors $u = [ 1 , 2 , 3 , 4 ]^\top$ and
$v = [ 5 , 6 , 7 , 8 ]^\top$, the concatenated vectors are
\begin{equation}
  [u, v]_0
  = 
  \begin{bmatrix} 1 \\ 5 \\ 2 \\ 6 \\ 3 \\ 7 \\ 4 \\ 8 \end{bmatrix},
  \quad
  [u, v]_1
  = 
  \begin{bmatrix} 1 \\ 2 \\ 5 \\ 3 \\ 6 \\ 4 \\ 7 \\ 8 \end{bmatrix},
  \quad
  [u, v]_2
  = 
  \begin{bmatrix} 1 \\ 2 \\ 3 \\ 5 \\ 4 \\ 6 \\ 7 \\ 8 \end{bmatrix}.
\end{equation}
The operation simply permutes the entries and places them into a vector, so it
is invertible.  See Fig.~\ref{fig:matrix_structure} for a graphical illustration. 

\begin{figure}
  \centering
  \begin{tabular}{|l|l|}
    \hline
    \begin{minipage}{0.7\textwidth}
    \includegraphics[width=1.0\textwidth]{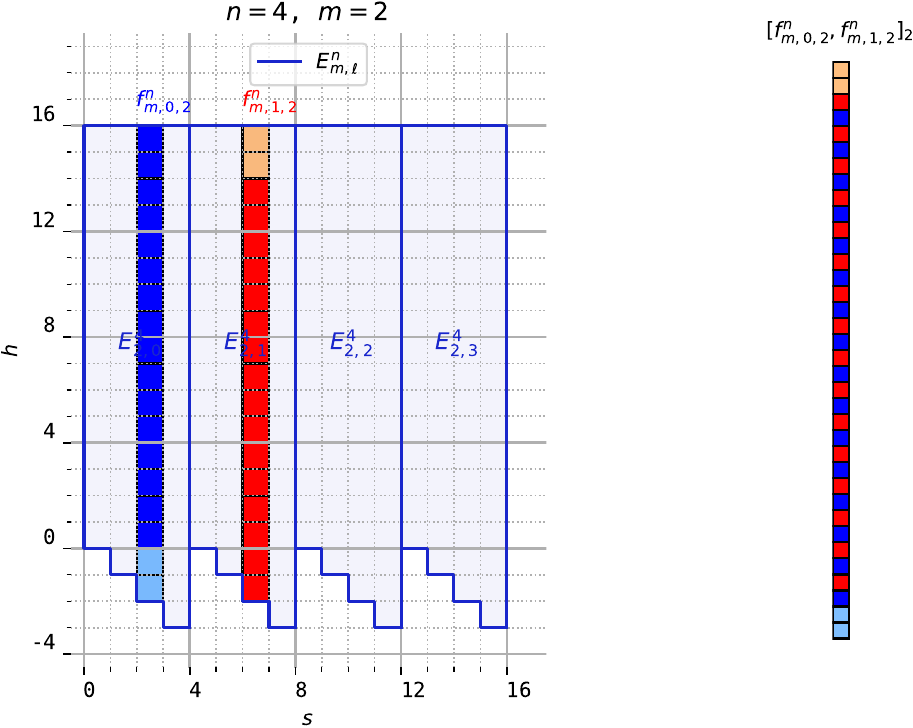}
    \end{minipage}
    \\
    \hline
    \begin{minipage}{0.7\textwidth}
    \includegraphics[width=1.0\textwidth]{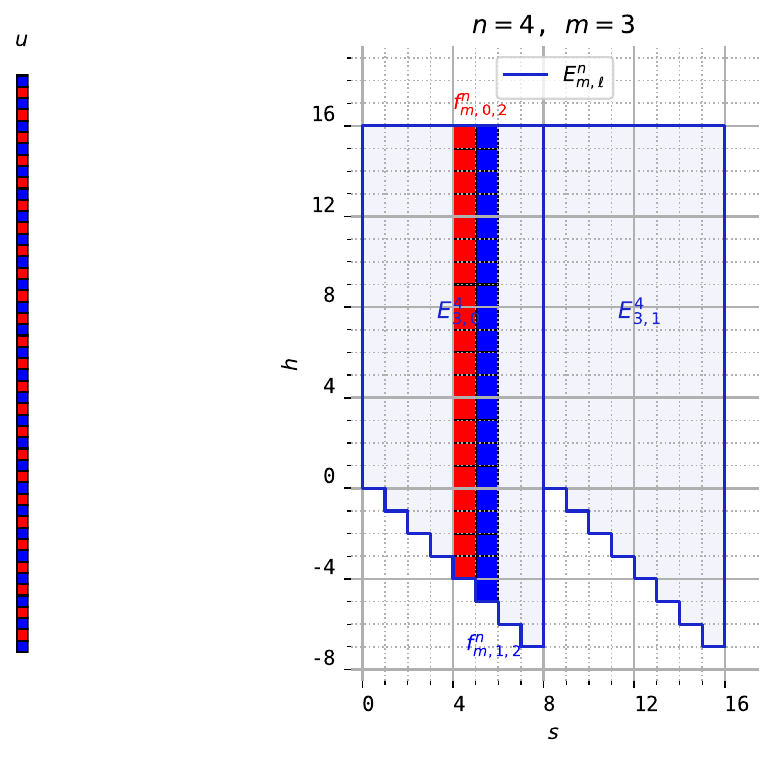}
    \end{minipage}
    \\
    \hline
  \end{tabular}
  \caption{Diagram depicting the flattening permutation $P^n_m$ given in
  \eqref{eq:P} (top), and diagram depicting the unflattening permutation $Q^n_m$
  given in \eqref{eq:Q} (bottom).}
  \label{fig:pq_permute}
\end{figure}

We will describe two permutations, the flattening permutation $P^n_m$
and the unflattening permutation $Q^n_m$, that will simplify
the data for the summation operations by $S^n_{(m)}$.
We define a flattening permutation to a
vector $P^n_m : \RR^{E^n_{m-1}} \to \RR^{\abs{E^n_{m-1}}}$. First, we define the
index
\begin{equation}
  \begin{aligned}
    k_{\ell, t}^{(P)}
    &:=
    2
    \dparen[\Big]{
    \sum_{\ell' < \ell} 
      |E^n_{m, \ell'}|
      +
    \sum_{t' < t} 
      \abs{E^n_{m,\ell, t'}}  
    },
    \qquad
    t \in [2^{m}].
  \end{aligned}
\end{equation} 
Then $P^n_m[f]$ for an image $f \in \RR^{E^n_{m-1}}$ is defined entrywise as
\begin{equation} \label{eq:P}
  \begin{aligned}
    \left( P^n_m [f] \right)_{j + k^{(P)}_{\ell, t}} 
    :=
    \left( \left[ f^n_{m, 2\ell, t} \mid f^n_{m, 2\ell + 1, t} \right]_t \right)_j,
    &&
    j \in [k^{(P)}_{\ell, t+1} - k^{(P)}_{\ell, t}],
    \,
    t \in [2^{n-m} - 1].
  \end{aligned}
\end{equation}
A diagram depicting the permutation $P^n_m$ is shown in Fig.~\ref{fig:pq_permute}.
Next, we define an unflattening permutation $Q^n_m : \RR^{|E^n_{m}|} \to
\RR^{E^n_{m}}$. Define the index
\begin{equation}
  k^{(Q)}_{\ell, t}
  :=
  \sum_{\ell' < \ell} | E^n_{m, \ell'} |
  +
  \sum_{t' < 2t} | E^n_{m, \ell, t'} |,
  \quad
  t \in [2^{m-1}].
\end{equation}
Let $Q^n_m [v] := g$ where the sections of the image $g$ are given in terms 
of $v$ by
\begin{equation} \label{eq:Q}
  \begin{aligned}
  \left[g^n_{m+1, \ell, 2t} \mid g^n_{m+1, \ell, 2t+1 } \right]
  :=
  \left[
    (v_{j + k^{(Q)}_{\ell, 2t}})_{j \in [k^{(Q)}_{m+1, \ell, t+1} - k^{(Q)}_{m+1,\ell, t}]}
  \right]^{-1}_0
  \end{aligned}
\end{equation}
for $t \in [2^{m} - 1]$.
A diagram depicting the permutation $Q^n_m$ appears in Fig.~\ref{fig:pq_permute}.

\begin{lemma}[Block convolution form] \label{lem:block_conv}
A factor of the single-quadrant ADRT $S^n_m$ in the case $m > 1$ can be written
\begin{equation} \label{eq:Snm_Znm}
  S^n_m
  =
  Q^n_m
  Z^n_m
  P^n_m,
\end{equation}
where $Z^n_m$ is given by \eqref{eq:Znm}, and the
permutation matrices $Q^n_m$ and $P^n_m$ are given by \eqref{eq:Q} and \eqref{eq:P}, respectively.
\end{lemma}

Now, we know the SVD of $Z^n_m$ \eqref{eq:Znm_svd} and how $Z^n_m$ is
related to $S^n_m$ \eqref{eq:Snm_Znm}, thus the SVD of the single-quadrant
factor $S^n_m$ follows immediately. We have that
\begin{equation} \label{eq:sq_svd}
  S^n_m
  =
  \left( Q^n_m U^n_m \right) 
  \Sigma^n_m 
  \left( P^{n \top}_m V^n_m \right)^\top,
  \quad
  m = 1, ...\,,n.
\end{equation}
It is easy to see that the SVD of the full ADRT factor $S^n_{(m)}$ follows too.

\begin{theorem}[SVD of $S^n_{(m)}$ for $m > 1$] 
  \label{prop:svd}
  The SVD of $S^n_{(m)}$ up to the ordering of the singular values is given by
  \begin{equation}
    S^n_{(m)}
    =
    U^n_{(m)}
    \Sigma^n_{(m)}
    V^{n \top}_{(m)}
  \end{equation}
  where
    $U^n_{(m)} := \Id_{4} \otimes \left( Q^n_m U^n_m \right),$
    $\Sigma^n_{(m)} := \Id_{4} \otimes \Sigma^n_m,$
    and
    $V^n_{(m)} := \Id_{4} \otimes \left( P^{n \top}_m V^n_m \right).$
  The multiplication by each factor can be computed in at most $\cO(N^2 \log N)$
  operations.
\end{theorem}
  
Since $U^n_{(m)}$ and $V^n_{(m)}$ are made up of block matrices $U^n_m$
and $V^n_m$, and these matrices are made up of columns from DST matrices,
the complexity of multiplying by these matrices (or their transposes) is
$\cO(N^2 \log N)$.

\section{Spectral decomposition of the cross-quadrant level}
\label{sec:spife}

Here we derive a spectral decomposition of the first factor $S^n_{(1)}$
\eqref{eq:Sn1} of the ADRT (Def.~\ref{def:Snm}).  The technique required for the first
level factor $S^n_{(1)}$ is different from that of the subsequent levels; we discuss the
difference and explain the difficulty. Regardless, we derive an alternative
spectral decomposition by employing a block-sinusoidal decomposition of the
original image.

\subsection{Cross-quadrant level vs. single-quadrant levels}

Recall that the first level factor $S^n_{(1)}$ of the ADRT flips and rotates the
original image via the permutations $T^n_\RN{1}, ...\,, T^n_\RN{4}$
\eqref{eq:permute}. We illustrate below that these operations prevent
$S^n_{(1)}$ from being a convolution-modulo-permutation. This fact implies that
the convolution theorem (Lem.~\ref{lem:conv}) cannot be applied to $S^n_{(1)}$,
unlike in the subsequent levels.

Let us compute the Gram matrix
\begin{equation}\label{eq:Sn1gram}
S^{n \top}_{(1)} S^n_{(1)} \in \RR^{N^2 \times N^2}
\quad \text{where} \quad
N = 2^n,
\end{equation}
and inspect the individual rows. We categorize each row depending on whether
its index pair (a pair in the 2D-image-indexing $[N/2] \times [N/2]$) has 
an
even ($e$) or odd ($o$) number as its first and second indices. Then each row
belongs to four possible categories depending on its index pair, $\{ (e, e), (e,
o), (o, e), (o, o) \}$. Reshaping the rows of each category (which are vectors
in $\RR^{N^2}$) into $N \times N$ images and displaying the non-zero entries
only, we obtain the $3 \times 3$ stencils shown in the first row of
Tab.~\ref{tab:stencil}. There are four different types of stencils, so there is no
clear convolutional structure.

Noting that the Gram matrix \eqref{eq:Sn1gram} can be expanded quadrant-wise
\begin{equation} \label{eq:SnTSn}
  S^{n \top}_{(1)} S^n_{(1)} 
  = 
  \sum_{j \in \{\RN{1}, ... ,
\RN{4}\}}
  (S^n_1 T^n_j )^\top (S^n_1 T^n_j),
\end{equation}
we also show in the subsequent rows in Tab.~\ref{tab:stencil} the contributions from
the individual terms, each coming from each of the four quadrants.
\begin{table}
  \centering
 \text{ 
\begin{tabular}{|r|c|c|c|c|}
  \hline
  &
  $(e, e)$
    &
  $(e, o)$
    &
  $(o, e)$
    &
  $(o, o)$
  \\
  \hline
  \hline
  Stencil
  &
  $\begin{matrix}
      1 & 2 & 2 \\
        & 8 & 2 \\
        &   & 1
    \end{matrix}$
    &
  $\begin{matrix}
        &   & 1 \\
        & 8 & 2 \\
      1 & 2 & 2
    \end{matrix}$
    &
  $\begin{matrix}
      2 & 2 & 1 \\
      2 & 8 &   \\
      1 &   &  
    \end{matrix}$
    &
  $\begin{matrix}
      1 &   &   \\
      2 & 8 &   \\
      2 & 2 & 1
    \end{matrix}$
  \\
  \hline
  \hline
  Quadrant $\RN{1}$
  &
  $\begin{matrix}
      1 & 1 &   \\
        & 2 &   \\
        &   &
    \end{matrix}$
    &
  $\begin{matrix}
        &   &   \\
        & 2 &   \\
        & 1 & 1
    \end{matrix}$
    &
  $\begin{matrix}
      1 & 1 &   \\
        & 2 &   \\
        &   &  
    \end{matrix}$
    &
  $\begin{matrix}
        &   &   \\
        & 2 &   \\
        & 1 & 1
    \end{matrix}$
  \\ \hline
  Quadrant $\RN{2}$
  &
  $\begin{matrix}
        &   &   \\
        & 2 & 1 \\
        &   & 1
    \end{matrix}$
    &
  $\begin{matrix}
        &   &   \\
        & 2 & 1 \\
        &   & 1 
    \end{matrix}$
    &
  $\begin{matrix}
      1 &   &   \\
      1 & 2 &   \\
        &   &
    \end{matrix}$
    &
  $\begin{matrix}
      1 &   &   \\
      1 & 2 &   \\
        &   &  
    \end{matrix}$
  \\ \hline
  Quadrant $\RN{3}$
  &
  $\begin{matrix}
        &   & 1 \\
        & 2 & 1 \\
        &   &
    \end{matrix}$
    &
  $\begin{matrix}
        &   & 1 \\
        & 2 & 1 \\
        &   &
    \end{matrix}$
    &
  $\begin{matrix}
        &   &   \\
      1 & 2 &   \\
      1 &   &   
    \end{matrix}$
    &
  $\begin{matrix}
        &   &   \\
      1 & 2 &   \\
      1 &   &  
    \end{matrix}$
  \\ \hline
  Quadrant $\RN{4}$
  &
  $\begin{matrix}
        & 1 & 1 \\
        & 2 &   \\
        &   &
    \end{matrix}$
    &
  $\begin{matrix}
        &   &   \\
        & 2 &   \\
      1 & 1 &    
    \end{matrix}$
    &
  $\begin{matrix}
        & 1 & 1 \\
        & 2 &   \\
        &   &   
    \end{matrix}$
    &
  $\begin{matrix}
        &   & \\
        & 2 & \\
      1 & 1 &
    \end{matrix}$
  \\
  \hline
 \end{tabular}
}
\caption{A description of the rows in $S^{n \top}_{(1)} S^{n}_{(1)}$. The first
row shows the ``stencil'' of each row as categorized by whether the placement
of the stencil has even or odd indexing in $[N/2] \times [N/2]$. The subsequent rows
detail the individual contributions coming from each of the four quadrants.
These stencils are different depending on the row category, so $S^n_{(1)}$ is
not a convolution.}
\label{tab:stencil}
\end{table}

Again, there are four different types of stencils depending on the quadrant and
the position of the row, so these individual operations are not convolutions.
Naturally, the convolution theorem does not apply in this situation.

We can describe this difference of the first level $S^n_{(1)}$ in another way.
Consider the na\"ive inverse that can be derived from the SVD of the
single-quadrant versions $S^n_1$ \eqref{eq:sq_Snm}.  Recall that the
single-quadrant ADRT at each level has the explicit SVD \eqref{eq:sq_svd}, e.g.,
for the first level we have $S^n_1 = U^n_1 \Sigma^n_1 V^{n\top}_1$.  Using this
directly in the definition of $S^n_{(1)}$ \eqref{eq:Sn1}, one derives a
straightforward decomposition
\begin{equation}
  \begin{aligned}
    S^n_{(1)}
    &=
    (\Id_4 \otimes U^n_1 )
    (\Id_4 \otimes \Sigma^n_1 )
    \left[
       T_\RN{1}^\top V^{n}_1 
       \mid
       T_\RN{2}^\top V^{n}_1
       \mid
       T_\RN{3}^\top V^{n}_1
       \mid
       T_\RN{4}^\top V^{n}_1 
    \right]^\top.
  \end{aligned}
\end{equation}
However, the last factor 
    $\left[
       T_\RN{1}^\top V^{n}_1 
       \mid
       T_\RN{2}^\top V^{n}_1
       \mid
       T_\RN{3}^\top V^{n}_1
       \mid
       T_\RN{4}^\top V^{n}_1 
    \right]$
has orthogonal rows but non-orthogonal columns when viewed as a matrix: so this
is not an SVD of $S^n_{(1)}$.  Still, one can use this decomposition to compute
the pseudo-inverse,
\begin{equation} \label{eq:spife-sq}
  S^{n \dagger}_{(1)}
  =
  \frac{1}{4}
  \left[
       T_\RN{1}
       \mid
       T_\RN{2}
       \mid
       T_\RN{3}
       \mid
       T_\RN{4}
    \right]^\top
    \left(\Id_4 \otimes V^n_1 \right)
    \left(\Id_4 \otimes (\Sigma^n_1)^{-1} \right)
    \left(\Id_4 \otimes U^{n\top}_1 \right),
\end{equation}
which amounts to computing the single-quadrant inverses for each quadrant
separately, undoing the permutations $T_\RN{1}, ...\, , T_\RN{4}$, then
averaging them pointwise. While this formula is algebraically correct, it has a
downside in terms of its stability: the diagonal matrix $\Id_4 \otimes
(\Sigma^n_1)^{-1}$ contains large numerical values, similar to the subsequent levels
$\Id_4 \otimes (\Sigma^n_m)^{-1}$. Therefore, the stability of this na\"ive
pseudo-inverse is equivalent to that of the exact inversion formula for the
single-quadrant ADRT \cite{OliviaGarcia2022}; asymptotically, it is as unstable
as the limited angle problem.

\subsection{Spectral decomposition of $S^n_{(1)}$}

Due to the observation laid out in the previous section, a different treatment
is needed to derive a decomposition for $S^n_{(1)}$. Instead of relying on the
convolution theorem (Lem.~\ref{lem:block_conv}), we derive an alternative explicit
decomposition. The key idea is to construct a 2D version of the DST basis that
is sinusoidal across $2 \times 2$ blocks rather than over individual entries.

\begin{figure}
  \centering
  \begin{tabular}{l}
  (a)
  \\
  \includegraphics[width=0.65\textwidth]{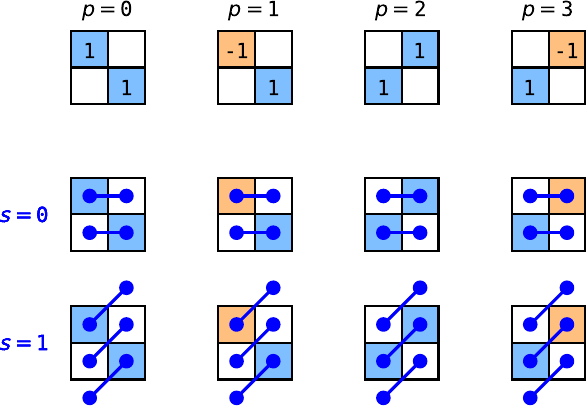}
  \\
  \\
  (b)
  \\
  \includegraphics[width=0.9\textwidth]{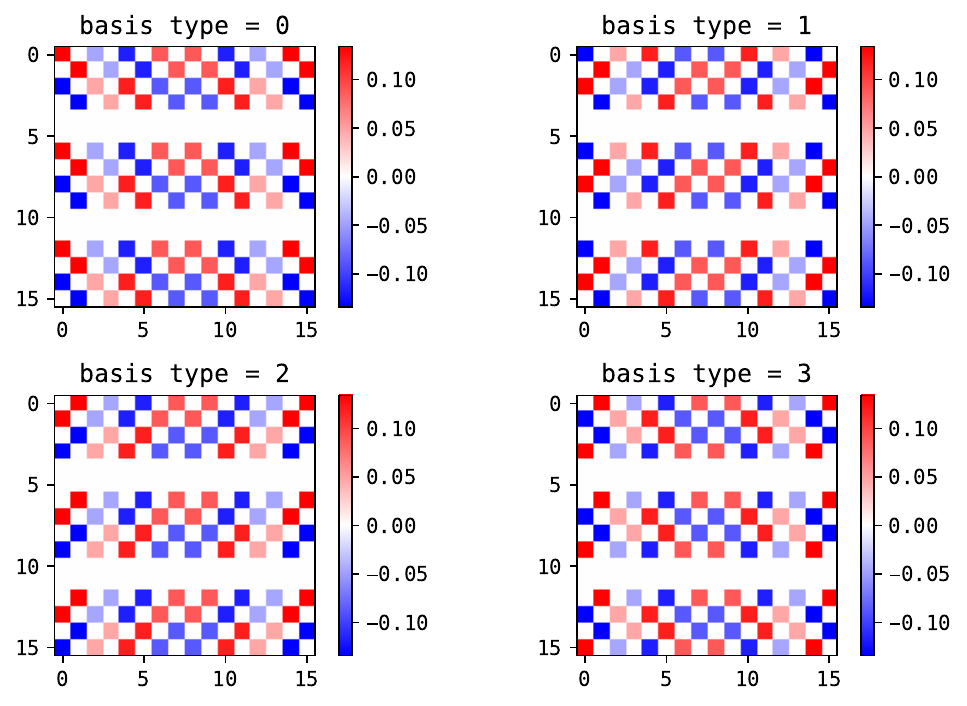}
  \end{tabular}
  \caption{(a) The $2 \times 2$ blocks of the basis functions for each type $p \in [4]$, and (b) an example of $v^{(p)}_{(k_1, k_2)}$ with $(k_1, k_2) = (5, 4)$ for each $p \in [4]$.}
  \label{fig:chi_block}
\end{figure}

We start by constructing a certain basis that will be useful. We decompose the
linear space of images $\RR^{N \times N}$ (where $N=2^n$) defined on the grid
$[N] \times [N]$  into four different types indexed by $p \in
[4]$, and each of these types into $N^2/4$ different sub-types indexed by $(j_1, j_2) \in [N/2] \times [N/2]$.
Define $\chi^{(p)}_{(j_1,j_2)} : [N] \times [N] \to \RR$ 
 over the grid $(i_1, i_2) \in [N] \times [N]$ with $p \in [4]$ and
$(j_1, j_2) \in [N/2] \times [N/2]$  as follows,
\begin{equation}
  \begin{aligned}
    &\chi^{(0)}_{(j_1, j_2)} (i_1, i_2)=
    &
    &\chi^{(1)}_{(j_1, j_2)} (i_1, i_2)=
    \\
     &
     \begin{cases}
      1 & (i_1, i_2) = (j_1, j_2),
      \\
      0 & (i_1,i_2) = (j_1, j_2 + 1),
      \\
      0 & (i_1, i_2) = (j_1+1, j_2),
      \\
      1 & (i_1, i_2) = (j_1+1, j_2 + 1),
      \\
      0 & \text{ otherwise},
     \end{cases} 
     &
     &
    \begin{cases}
    -1 & (i_1, i_2) = (j_1, j_2),
     \\  
     0 & (i_1,i_2) = (j_1, j_2 + 1),
     \\  
     0 & (i_1, i_2) = (j_1+1, j_2),
     \\  
     1 & (i_1, i_2) = (j_1+1, j_2 + 1),
     \\
     0 & \text{ otherwise},
    \end{cases} 
    \\
    &\chi^{(2)}_{(j_1, j_2)} (i_1, i_2) =
    &
    &\chi^{(3)}_{(j_1, j_2)} (i_1, i_2) =
    \\
   &
   \begin{cases}
    0 & (i_1, i_2) = (j_1, j_2),       
    \\  
    1 & (i_1,i_2) = (j_1, j_2 + 1),    
    \\  
    1 & (i_1, i_2) = (j_1+1, j_2),     
    \\  
    0 & (i_1, i_2) = (j_1+1, j_2 + 1), 
    \\
    0 & \text{ otherwise},
   \end{cases} 
   &
   &
    \begin{cases}
    0 & (i_1, i_2) = (j_1, j_2), 
    \\  
    -1 & (i_1,i_2) = (j_1, j_2 + 1),    
    \\  
    1 & (i_1, i_2) = (j_1+1, j_2),     
    \\  
    0 & (i_1, i_2) = (j_1+1, j_2 + 1), 
    \\
    0 & \text{ otherwise}.
    \end{cases}
  \end{aligned}
\end{equation}
Each of these basis functions is supported in $2\times2$ blocks.  An
illustration of the blocks is shown in Fig.~\ref{fig:chi_block}.
Next, we define the set of functions
$\chi^{(p)} \subset (\RR^{N \times N})^{(N/2) \times (N/2)}$  
based on their type $p$,
\begin{equation}
  \chi^{(p)}
  :=
  \left\{
      \chi^{(p)}_{(j_1, j_2)}
      \, \mid \,
    (j_1, j_2) \in [N/2] \times [N/2]
  \right\},
  \quad
  p \in [4].
\end{equation}

Then clearly, the union of their span is the linear space of $N \times N$
images,
  $\Span \bigcup_{p \in [4]} \chi^{(p)}
  =
  \RR^{N \times N},$
and these functions are pairwise orthogonal in the sense that
\begin{equation}
  \sum_{(i_1, i_2) \in [N] \times [N]}
  \chi^{(p)}_{j_1, j_2} (i_1, i_2)
  \chi^{(p')}_{j_1', j_2'} (i_1, i_2)
  =
  2 \delta_{pp'} \delta_{j_1 j_1'} \delta_{j_2 j_2'},
\end{equation}
for $p, p' \in [4]$, $j_1, j_1', j_2, j_2' \in [N/2]$. Hence, $\cup_{p \in [4]}
\chi^{(p)}$ is an orthogonal basis for $\RR^{N \times N}$.
Next, we define certain sinusoidal functions. Let us define for each $(k_1, k_2)
\in [N/2] \times [N/2]$ the following basis functions for $\RR^{N/2 \times N/2}$,
\begin{equation}
  \begin{aligned}
  \phi_{(k_1, k_2)} (j_1, j_2)
  &:=
  \frac{1}{N + 2}
  \sin \left( \frac{k_1 \pi j_1}{(N/2) + 1} \right)
  \sin \left( \frac{k_2 \pi j_2}{(N/2) + 1} \right),
  \end{aligned}
\end{equation}
and their variants
\begin{equation} \label{eq:psi}
  \begin{aligned}
  \psi^{+ \circ}_{(k_1, k_2)} (j_1, j_2)
  &:=
  \frac{c_{k_1}}{\sqrt{N + 2}}
  \sin \left( \frac{k_1 \pi (j_1 + \half)}{(N/2) + 1} \right)
  \sin \left( \frac{k_2 \pi j_2}{(N/2) + 1} \right),
  \\
  \psi^{\circ +}_{(k_1, k_2)} (j_1, j_2)
  &:=
  \frac{c_{k_2}}{\sqrt{N + 2}}
  \sin \left( \frac{k_1 \pi j_1}{(N/2) + 1} \right)
  \sin \left( \frac{k_2 \pi (j_2 + \half)}{(N/2) + 1} \right),
  \\
  \psi^{- \circ}_{(k_1, k_2)} (j_1, j_2)
  &:=
  \frac{c_{k_1}}{\sqrt{N + 2}}
  \cos \left( \frac{k_1 \pi (j_1 + \half)}{(N/2) + 1} \right)
  \sin \left( \frac{k_2 \pi j_2}{(N/2) + 1} \right),
  \\
  \psi^{\circ -}_{(k_1, k_2)} (j_1, j_2)
  &:=
  \frac{c_{k_2}}{\sqrt{N + 2}}
  \sin \left( \frac{k_1 \pi j_1}{(N/2) + 1} \right)
  \cos \left( \frac{k_2 \pi (j_2 + \half)}{(N/2) + 1} \right),
  \end{aligned}
\end{equation}
in which $c_k$ is the normalization factor for DCT-2 and DST-2, 
  $c_k 
  =
  \begin{cases}
    \frac{1} {\sqrt{N}} & \text{ if } k = 0,\\
    \frac{1} {\sqrt{2N}} & \text{ otherwise}.\\
  \end{cases}$

Finally, define the basis vectors
$v^{(p)}_{(k_1,k_2)} : [N] \times [N] \to \RR$ where $p \in [4]$ and
$(k_1, k_2) \in [N/2] \times [N/2]$ given by
\begin{equation} \label{eq:vpk1k2}
    v^{(p)}_{(k_1, k_2)}
    =
    \frac{1}{\sqrt{2}}
    \sum_{(j_1, j_2) \in [N/2] \times [N/2]}
      \phi_{(k_1, k_2)} (j_1, j_2) \, \chi^{(p)}_{(j_1, j_2)}.
\end{equation}
These functions are pairwise orthogonal in the sense that
\begin{equation}
  \sum_{(i_1, i_2) \in [N] \times [N]}
  v^{(p)}_{(k_1, k_2)} (i_1, i_2)
  v^{(p')}_{(k_1', k_2')} (i_1, i_2)
  =
  \delta_{pp'} \delta_{k_1 k_1'} \delta_{k_2 k_2'},
\end{equation}
for $p, p' \in [4]$, $k_1, k_1', k_2, k_2' \in [N/2]$. 
This is due to the orthogonality of the sets of functions 
$( \phi_{(k_1,k_2)} )_{(k_1, k_2) \in [N/2] \times [N/2]}$
and 
$( \chi^{(p)}_{(j_1,j_2)} )_{p \in [4], (j_1, j_2) \in [N/2] \times [N/2]}$
that yields
\begin{equation}
  \begin{aligned}
  &\sum_{(i_1, i_2) \in [N] \times [N]}
  v^{(p)}_{(k_1, k_2)} (i_1, i_2)
  v^{(p')}_{(k_1', k_2')} (i_1, i_2)
  \\ &=
  \frac{1}{2}
  \sum_{(i_1, i_2)}
  \sum_{(j_1, j_2) \in [N/2] \times [N/2]}
  \phi_{(k_1, k_2)}   (j_1,j_2)\chi^{(p)}_{(j_1, j_2)} (i_1, i_2)
  \phi_{(k_1', k_2')} (j_1,j_2)\chi^{(p')}_{(j_1, j_2)} (i_1, i_2)
  \\ &=
  \frac{1}{2}
  \sum_{(j_1, j_2)}
  \phi_{(k_1, k_2)}   (j_1,j_2)
  \phi_{(k_1', k_2')} (j_1,j_2)
  \left[
  \sum_{(i_1, i_2)}
  \chi^{(p)}_{(j_1, j_2)} (i_1, i_2)
  \chi^{(p')}_{(j_1, j_2)} (i_1, i_2)
  \right]
  \\ &=
  \frac{1}{2}
  \sum_{(j_1, j_2)}
  \phi_{(k_1, k_2)}   (j_1,j_2)
  \phi_{(k_1', k_2')} (j_1,j_2)
  (2\delta_{pp'})
  \\ &=
  \delta_{pp'} \delta_{k_1 k'_1} \delta_{k_2 k'_2}
  \end{aligned}
\end{equation}

This basis will be used
for our spectral decomposition.
Let $O^n_1 $ be a permutation of images in $\RR^{\abs{E^n_1}}$ reordering each
column into four blocks,
\begin{equation} \label{eq:ON}
  O^n_1 v
  :=
  \begin{bmatrix}
    \mu_1
    \\
    \mu_2
    \\
    \mu_3
    \\
    \mu_4
  \end{bmatrix},
  \qquad
  \mu_1, \mu_2, \mu_3 \in \RR^{N/2 \times N/2},
  \quad
  \mu_4 \in \RR^{((N/2) + 1) \times N/2},
\end{equation}
where the permuted entries are given by
\begin{equation} 
  \begin{aligned}
  &(\mu_i)(j_1,j_2) := v \dparen[\big]{(i \brem 2) + 2(j_1-1),\, \lfloor i / 2 \rfloor + j_2},
  \\
  &\qquad \text{for } i \in [4]_1,
  \quad
  j_1 \in 
  \left.
  \begin{cases}
    [N/2] & \text{ if } i < 4\\
    [(N/2) +1 ] & \text{ if } i =4\\
  \end{cases}
  \right\},
  \quad
  j_2 \in [N/2].
  \end{aligned}
\end{equation}
Since it is a permutation, its inverse is equal to its transpose,
$(O^n_1)^{-1} = O^{n \top}_1$. We use this transpose in the following.

\begin{theorem}[Spectral decomposition of $\npmzS{1}$] \label{thm:spd}
We have the decomposition
\begin{equation}
  S^n_{(1)}
  =
  U_{(1)}^n \Sigma_{(1)}^n V_{(1)}^{n \top},
\end{equation}
where the columns of $U_{(1)}^n$ and $V_{(1)}^n$ will be indexed by
$(p, k_1, k_2) \in [4]\times[N/2]\times[N/2]$ and are given as follows:
The columns of $V_{(1)}^n$ are $v_{(k_1, k_2)}^{(p)}$ \eqref{eq:vpk1k2}
and the columns of $U_{(1)}^n$ are
\begin{equation}
  u^{(p)}_{(k_1, k_2)}
  :=
  \begin{bmatrix}
  O^{n \top}_1 &              &              &              \\
               & O^{n \top}_1 &              &              \\
               &              & O^{n \top}_1 &              \\
               &              &              & O^{n \top}_1 \\
  \end{bmatrix}
  w^{(p)}_{(k_1, k_2)}
\end{equation}
where 
\begin{equation}\label{eq:ulvl1}
  \begin{aligned}
  \begin{matrix}
  w^{(0)}_{(k_1, k_2)}=
  \\
  \frac{1}{\sqrt{2}}
  \begin{bcolarray}
    \phi_{(k_1, k_2)}
    \\
    \phi_{(k_1, k_2)}
    \\
    2\phi_{(k_1, k_2)}
    \\
    \\
    \hline
    \phi_{(k_1, k_2)}
    \\
    \phi_{(k_1, k_2)}
    \\
    2\phi_{(k_1, k_2)}
    \\
    \\
    \hline
    \phi_{(k_1, k_2)}
    \\
    \phi_{(k_1, k_2)}
    \\
    
    \\
    \sigma^+_{k_1}
    \psi^{+ \circ}_{(k_1, k_2)}
    \\
    \hline
    \phi_{(k_1, k_2)}
    \\
    \phi_{(k_1, k_2)}
    \\
    
    \\
    \sigma^+_{k_2} 
    \psi^{\circ +}_{(k_1, k_2)}
  \end{bcolarray}
  \end{matrix},
  \begin{matrix}
  w^{(1)}_{(k_1, k_2)} =
  \\
  \frac{1}{\sqrt{2}}
  \begin{bcolarray}
    \phi_{(k_1, k_2)}
    \\
    -\phi_{(k_1, k_2)}
    \\
    
    \\
    
    \\
    \hline
    \phi_{(k_1, k_2)}
    \\
    -\phi_{(k_1, k_2)}
    \\
    
    \\
    
    \\
    \hline
    -\phi_{(k_1, k_2)}
    \\
    \phi_{(k_1, k_2)}
    \\
    
    \\
    -\sigma^-_{k_1} 
    \psi^{- \circ}_{(k_1, k_2)}
    \\
    \hline
    \phi_{(k_1, k_2)}
    \\
    -\phi_{(k_1, k_2)}
    \\
    
    \\
    -\sigma^-_{k_2} 
    \psi^{\circ -}_{(k_1, k_2)}
  \end{bcolarray}\end{matrix},
  \,
  \begin{matrix}
  w^{(2)}_{(k_1, k_2)} =
  \\
  \frac{1}{\sqrt{2}}
  \begin{bcolarray}
    \phi_{(k_1, k_2)}
    \\
    \phi_{(k_1, k_2)}
    \\
    
    \\
    \sigma^+_{k_2}
    \psi^{\circ +}_{(k_1, k_2)}
    \\
    \hline
    \phi_{(k_1, k_2)}
    \\
    \phi_{(k_1, k_2)}
    \\
    
    \\
    \sigma^+_{k_1}
    \psi^{+ \circ}_{(k_1, k_2)}
    \\
    \hline
    \phi_{(k_1, k_2)}
    \\
    \phi_{(k_1, k_2)}
    \\
    2\phi_{(k_1, k_2)}
    \\
    
    \\
    \hline
    \phi_{(k_1, k_2)}
    \\
    \phi_{(k_1, k_2)}
    \\
    2\phi_{(k_1, k_2)}
    \\
    
    \\
  \end{bcolarray} \end{matrix},
  \begin{matrix} w^{(3)}_{(k_1, k_2)} =
    \\
    \frac{1}{\sqrt{2}}
  \begin{bcolarray}
    \phi_{(k_1, k_2)}
    \\
    -\phi_{(k_1, k_2)}
    \\
    
    \\
    -\sigma^-_{k_2}
    \psi^{\circ -}_{(k_1, k_2)}
    \\
    \hline
    -\phi_{(k_1, k_2)}
    \\
    \phi_{(k_1, k_2)}
    \\
    \\
    \sigma^-_{k_1}
    \psi^{- \circ}_{(k_1, k_2)}
    \\
    \hline
    \phi_{(k_1, k_2)}
    \\
    -\phi_{(k_1, k_2)}
    \\
    \\
    \\
    \hline
    \phi_{(k_1, k_2)}
    \\
    -\phi_{(k_1, k_2)}
    \\
    
    \\
    
    \\
  \end{bcolarray} \end{matrix}.
  \end{aligned} 
\end{equation} 
The diagonal values of $\Sigma^n_{(1)}$ are given by
\begin{equation} \label{eq:slvl1}
  \begin{aligned}
    \sigma_{(k_1,k_2)}^{(0)}
    =
    \sigma_{(k_1,k_2)}^{(2)}
    &=
    \dparen[\big]{ 8 + \half (\sigma_{k_1}^+)^2 + \half (\sigma_{k_2}^+)^2}^\half,
    \\
    \sigma_{(k_1,k_2)}^{(1)}
    =
    \sigma_{(k_1,k_2)}^{(3)}
    &=
    \dparen[\big]{4 + \half (\sigma_{k_1}^-)^2 + \half (\sigma_{k_2}^-)^2 }^\half,
  \end{aligned}
\end{equation}
where
  $\sigma^+_k := 2 \cos \left( \frac{k \pi}{N + 2} \right)$
  and
  $\sigma^-_k := 2 \sin \left( \frac{k \pi}{N + 2} \right).$
\end{theorem}

\begin{proof}

It suffices to consider the case for one basis type $p$ and one quadrant, as the
computation for the other cases is similar.  We choose basis type $p=3$ and
quadrant $\RN{1}$. 

Fix $(k_1, k_2) \in [N/2] \times [N/2]$. Consider the permuted vector
\begin{equation}
  T_\RN{1} v^{(p)}_{(k_1, k_2)}.
\end{equation}
Viewed as a 2D image, the $2 \times 2$ subblocks are now scalar multiples of the
skew symmetric matrix
  $\begin{bmatrix}0 & -1 \\ 1 & 0\end{bmatrix}.$
The calculation of the restriction of $u^{(p)}_{(k_1, k_2)}$ to 
quadrant~\RN{1}, which is simply
  $\left. u^{(p)}_{(k_1, k_2)} \right\rvert_{\RN{1}} (h, s)
  =
  [S^n_{1} T_{\RN{1}} v^{(p)}_{(k_1,k_2)}] (h, s),$
  for
  $(h, s) \in E^n_1$,
can be decomposed into four parts, depending on whether $h$ and $s$ are each
even or odd. We consider each case:
When $s$ is even, the sum occurs only within each subblock, and the sum is
trivial since it adds zero to one entry in the subblock. So, when $h$ is even,
one obtains
\begin{equation}
  u(h, s)
  =
- \frac{1}{\sqrt{2}} \phi_{(k_1, k_2)} 
  \dparen[\big]{\lfloor s/2 \rfloor, (N/2) - 1 - \lfloor h/2 \rfloor},
\end{equation}
and when $h$ is odd,
\begin{equation}
  u(h, s)
  =
  \frac{1}{\sqrt{2}} \phi_{(k_1, k_2)}
  \dparen[\big]{\lfloor s/2 \rfloor, (N/2) - 1 - \lfloor h/2 \rfloor}.
\end{equation}
Next, when $s$ is odd and $h$ is odd, we are only summing zero entries, so
\begin{equation}
  u(h, s)
  =
  0,
\end{equation}
and finally when $h$ is even,
\begin{equation}
  \begin{aligned}
    u(h, s)
    &=
    \phi_{(k_1, k_2)} 
   \dparen[\big]{\lfloor s/2 \rfloor, (N/2) - \lfloor h/2 \rfloor}
   \\
   & \qquad
     -
    \phi_{(k_1, k_2)} 
    \dparen[\big]{\lfloor s/2 \rfloor, (N/2) - 1 - \lfloor h/2 \rfloor}
    \\
    &=
    - \sigma_{k_2}^-
    \psi_{(k_1, k_2)}^{\circ -} 
    \dparen[\big]{\lfloor s/2 \rfloor, N/2 - \lfloor h/2 \rfloor},
  \end{aligned}
\end{equation}
by Lem.~\ref{lem:conv}. This calculation corresponds to the calculation of the top four
entries of $w^{(3)}_{(k_1,k_2)}$ in \eqref{eq:ulvl1} (see Fig.~\ref{fig:chi_block}
for a visual illustration).
The diagonal values \eqref{eq:slvl1} are obtained by taking the $\ell^2$-norm of
the vectors in \eqref{eq:ulvl1}.

\end{proof}

\section{Spectral pseudo-inverse of $S^n_{(1)}$}
\label{sec:spife_def}

We finally turn to the problem of solving the system
  $S^{n}_{(1)} [x] = b,$
  with $b \in \text{range}(S^{n}_{(1)}).$
Plugging in our spectral decomposition from Thm.~\ref{thm:spd}, we have
\begin{equation} \label{eq:S1svdsolve}
  U^{n}_{(1)} \Sigma^n_{(1)} V^{n \top}_{(1)} x = b.
\end{equation}
If this decomposition were the SVD, we would simply multiply by $U^{n
\top}_{(1)}$ on both sides. Unfortunately, despite the notation resembling an
SVD, the matrix $U^n_{(1)}$ here is not orthogonal in the standard inner
product. It is not difficult to compute the Gram matrix associated with
$U^n_{(1)}$ and see that it has the block matrix form
\begin{equation}
  U^{n \top}_{(1)} U^n_{(1)}
  =
  \begin{bmatrix}
  \Id_{N^2 + N/2} & B_{12} & B_{13}      &        \\    
  B_{12}^\top & \Id_{N^2 + N/2}    &             &        \\    
  B_{13}^\top &        & \Id_{N^2 + N/2}         & B_{34} \\    
              &        & B_{34}^\top & \Id_{N^2 + N/2}    \\    
  \end{bmatrix}
\end{equation}
where $B_{13}$ is diagonal, and the $k$-th diagonal of $B_{12}$ is zero for $k$
even. The blocks $B_{12}$ and $B_{34}$ only involve inner products between
$\psi^{\pm \circ}_{(k_1, k_2)}$ and $\psi^{\circ \pm}_{(k_1, k_2)}$, and the
diagonals of $B_{12}$ are $|\phi^{(p)}_{(k_1, k_2)}|^2$ up to constant
multiples.  So the columns of $U$ are orthogonal with respect to a weighted
inner product, with the weight $(U_{(1)}^\top U_{(1)})^{-\half}$, 
and this weight is not
diagonal. Hence, it might appear that more work would be necessary to invert this system.
However, the solution can still be computed efficiently.

\begin{theorem}[Explicit spectral pseudo-inverse of $\npmzS{1}$] \label{prop:ispife}
Given $b \in \rg \npmzS{1}$, the solution $x \in \RR^{N^2}$ to the linear problem $S^{n}_{(1)} [x] = b$
can be computed by applying two orthogonal transformations at the computational
complexity $\cO(N^2 \log N)$. More precisely, let $U_{(1)}^n$, $V^n_{(1)}$ be as
in Thm.~\ref{thm:spd}; then the solution is given by
\begin{equation} \label{eq:spife1}
  x = \widehat{V}^n_{(1)} \widehat{U}_{(1)}^{n \top} \hat{b},
\end{equation}
where the matrices on the RHS are 
  $\widehat{U}_{(1)}^n := M^n_{(1)} U^n_{(1)} \Sigma^n_{(1)} Q^n_{(1)},$
  $\widehat{V}_{(1)}^n := V^n_{(1)} Q^n_{(1)},$
  and
  $\hat{b} := M^n_{(1)} b.$
$M_{(1)}^n$ is a square diagonal matrix with 1s and 0s as entries (a
truncation matrix), and $Q^n_{(1)}$ is an orthogonal matrix.
\end{theorem}

\begin{proof}
First, we derive the orthogonal matrix $\widehat{U}^n_{(1)}$ through a
transformation of $U^n_{(1)}$. Denote by $\bar{w}^{(p)}_{(k_1, k_2)}$ the
restriction of $w^{(p)}_{(k_1, k_2)}$ \eqref{eq:ulvl1} where row-blocks involving
the variants \eqref{eq:psi} are deleted. That is,
\begin{equation}\label{eq:barulvl1}
  \begin{aligned}
  \begin{matrix}
  \bar{w}^{(0)}_{(k_1, k_2)}=
  \\
  \frac{1}{\sqrt{2}}
  \begin{bcolarray}
    \phi_{(k_1, k_2)}
    \\
    \phi_{(k_1, k_2)}
    \\
    2\phi_{(k_1, k_2)}
    \\ 
    \hline
    \phi_{(k_1, k_2)}
    \\
    \phi_{(k_1, k_2)}
    \\
    2\phi_{(k_1, k_2)}
    \\
    \hline
    \phi_{(k_1, k_2)}
    \\
    \phi_{(k_1, k_2)}
    \\
    \\
    \hline
    \phi_{(k_1, k_2)}
    \\
    \phi_{(k_1, k_2)}
    \\
    \\
  \end{bcolarray}
  \end{matrix},
  \quad
  \begin{matrix}
  \bar{w}^{(1)}_{(k_1, k_2)} =
  \\
  \frac{1}{\sqrt{2}}
  \begin{bcolarray}
    \phi_{(k_1, k_2)}
    \\
    -\phi_{(k_1, k_2)}
    \\
    
    \\
    \hline
    \phi_{(k_1, k_2)}
    \\
    -\phi_{(k_1, k_2)}
    \\
    
    \\
    \hline
    -\phi_{(k_1, k_2)}
    \\
    \phi_{(k_1, k_2)}
    \\
    
    \\
    \hline
    \phi_{(k_1, k_2)}
    \\
    -\phi_{(k_1, k_2)}
    \\
    
    \\
  \end{bcolarray} \end{matrix},
  \quad
  \begin{matrix}
  \bar{w}^{(2)}_{(k_1, k_2)} =
  \\
  \frac{1}{\sqrt{2}}
  \begin{bcolarray}
    \phi_{(k_1, k_2)}
    \\
    \phi_{(k_1, k_2)}
    \\
    
    \\
    \hline
    \phi_{(k_1, k_2)}
    \\
    \phi_{(k_1, k_2)}
    \\
    
    \\
    \hline
    \phi_{(k_1, k_2)}
    \\
    \phi_{(k_1, k_2)}
    \\
    2\phi_{(k_1, k_2)}
    \\
    \hline
    \phi_{(k_1, k_2)}
    \\
    \phi_{(k_1, k_2)}
    \\
    2\phi_{(k_1, k_2)}
    \\
  \end{bcolarray} \end{matrix},
  \quad
  \begin{matrix} \bar{w}^{(3)}_{(k_1, k_2)} =
    \\
  \frac{1}{\sqrt{2}}
  \begin{bcolarray}
    \phi_{(k_1, k_2)}
    \\
    -\phi_{(k_1, k_2)}
    \\
    
    \\
    \hline
    -\phi_{(k_1, k_2)}
    \\
    \phi_{(k_1, k_2)}
    \\
    \\
    \hline
    \phi_{(k_1, k_2)}
    \\
    -\phi_{(k_1, k_2)}
    \\
    \\
    \hline
    \phi_{(k_1, k_2)}
    \\
    -\phi_{(k_1, k_2)}
    \\
    
    \\
  \end{bcolarray} \end{matrix}.
  \end{aligned} 
\end{equation} 
Then we take unit vectors in the directions
\begin{equation} \label{eq:lincomb}
  \begin{aligned}
    \hat{w}^{(0)} &\parallel \bar{w}^{(0)} + \bar{w}^{(2)},
    &&
    % \dnorm[\big]{\hat{w}^{(0)}_{(k_1, k_2)}}_{2} = 1,
    \hat{w}^{(1)} &\parallel \bar{w}^{(1)},
    &&
    \hat{w}^{(2)} &\parallel \bar{w}^{(0)} - \bar{w}^{(2)},
    &&
    \hat{w}^{(3)} &\parallel \bar{w}^{(3)},
  \end{aligned}
\end{equation}
where the notation $\cdot \parallel \cdot$ is used to denote that the two
vectors are parallel, and we have omitted the subscripts. We obtain the
vectors 
\begin{align}
  \hat{w}^{(0)} \parallel
   \begin{bcolarray}
     \phi
      \\
     \phi
      \\
     \phi
      \\
      \hline
     \phi
      \\
     \phi
      \\
     \phi
      \\
      \hline
     \phi
      \\
     \phi
      \\
     \phi
      \\
      \hline
     \phi
      \\
     \phi
      \\
     \phi
      \\
    \end{bcolarray},
    &&
    \hat{w}^{(1)} \parallel
    \begin{bcolarray}
      \phi
      \\
     -\phi
      \\
      \\
      \hline
      \phi
      \\
     -\phi
      \\
      \\
      \hline
     -\phi
      \\
      \phi
      \\
      \\
      \hline
      \phi
      \\
     -\phi
      \\
      \\  
    \end{bcolarray},
    &&
    \hat{w}^{(2)}_{(k_1, k_2)} \parallel 
    \begin{bcolarray}
      \\
      \\
     \phi
      \\
      \hline
      \\
      \\
     \phi
      \\
      \hline
      \\
      \\
    -\phi
      \\
      \hline
      \\
      \\
    -\phi
      \\
    \end{bcolarray},
    &&
    \hat{w}^{(3)} \parallel 
    \begin{bcolarray}
      \phi
      \\
     -\phi
      \\
      \\
      \hline
     -\phi
      \\
      \phi
      \\
      \\
      \hline
      \phi
      \\
     -\phi
      \\
      \\
      \hline
      \phi
      \\
     -\phi
      \\
      \\
    \end{bcolarray}.
\end{align}

The linear combinations in \eqref{eq:lincomb} are performed by the matrix $Q^n_{(1)}$,
and the normalization is achieved by a rescaling using the matrix
$\Sigma^{n}_{(1)}$. That $\{\hat{w}^{(p)}_{(k_1, k_2)}\}_{p, k_1, k_2}$ forms an
orthonormal set of vectors follows from direct calculation.
So, our original linear system \eqref{eq:S1svdsolve} can be truncated as
$M^n_{(1)} U^n_{(1)} \Sigma^n_{(1)} V^{n \top}_{(1)} x = M^n_{(1)} b.$
Then, we rewrite
\begin{equation}
\underbrace{M^n_{(1)} U^n_{(1)} \Sigma^n_{(1)} Q^n_{(1)}}_{\widehat{U}^n_{(1)}}
\underbrace{Q^{n \top}_{(1)} V^{n \top}_{(1)}}_{\widehat{V}^{n \top}_{(1)}}
  x = \hat{b},
\end{equation}
and the orthogonal matrices on the LHS can be inverted, leading to
\eqref{eq:spife1}. Since it is made up of blocks of DSTs, this computation can be
done with complexity $\cO(N^2 \log N)$.
\end{proof}

Finally, using the SVD from Sec.~\ref{sec:svd} and the spectral decomposition
from Sec.~\ref{sec:spife}, we obtain the desired pseudo-inverse of the ADRT.  

\begin{theorem}[Spectral Pseudo-Inverse, Fast and Explicit (SPIFE)]
  \label{prop:spifefull}
Let $\widehat{U}^n_{(1)}$ be as in
Thm.~\ref{prop:ispife}, $V^n_{(1)}$ as in Thm.~\ref{thm:spd}, and $U_{(m)}^n,
\Sigma_{(m)}^n, V_{(m)}^n$ as in Thm.~\ref{prop:svd}. Then a pseudo-inverse
of $\nR$ is given by
\begin{equation} \label{eq:spife}
   (R^n)^\dagger
   =
   \widehat{V}^n_{(1)}
   \widehat{U}^{n \top}_{(1)}
   \prod_{m=2}^n
   V^n_{(m)}
   (\Sigma^n_{(m)})^{-1}
   U^{n \top}_{(m)}.
\end{equation}
The pseudo-inverse can be multiplied explicitly in $\cO(N^2 \log^2 N)$
operations.
\end{theorem}

\section{Numerical experiments} \label{sec:numerics}

Here we conduct numerical experiments regarding SPIFE.
We start by comparing the inversion results of five different methods:
\begin{description}[font={\ttfamily},leftmargin=\parindent+50pt,style=nextline,itemsep=2pt,labelindent=\parindent,rightmargin=30pt]
\item[spife] The pseudo-inverse using the spectral decomposition
  \eqref{eq:spife}
  \item[spife-sq] The pseudo-inverse using the single-quadrant SVDs and
  taking the average across four quadrants \eqref{eq:spife-sq}
  \item[cg] The conjugate gradient (CG) iteration
  \item[fmg] The full multigrid (FMG) method
  \item[alg] The algebraically exact inversion formula
\end{description}
The implementation of the experiments was done using the open source Python
software package {\tt adrt} \cite{Otness2023}.

For useful comparisons, we run both the CG and FMG only up to $\log N$
iterations, so that the computational cost is comparable; both iterations cost
at least $\sim N^2 \log N$ operations per iteration, so $\log N$ iterations
result in a total cost of the order $\cO(N^2 \log^2 N)$, the cost of our
pseudo-inverse. That is, in all our comparisons involving iterative methods
below, $\log N$ iterations were used for data corresponding to size $N \times
N$ (recall $\log$ here denotes the logarithm of base two). For example, for data from
an image of size $N=128$, $7$ iterations were used.
Note, however, that {\tt alg} has the slightly lower complexity
$\cO(N^2 \log N)$.
For a sense of benchmark performance for the iterative methods, we
refer the reader to the observation regarding error behavior for FMG from the work
\cite{Press06drt}. A steady reduction of error per iteration was shown therein,
following the pattern
$\Delta \log \text{error} (\text{per iteration}) \approx - 13.8/\log^2 N$ (Eq.
14 in \cite{Press06drt}).

\begin{figure}
  \centering
  \begin{tabular}{c|c}
    \hline
    \\
    (a)
    &
    \begin{minipage}{0.6\textwidth}
    \includegraphics[width=1.0\textwidth]{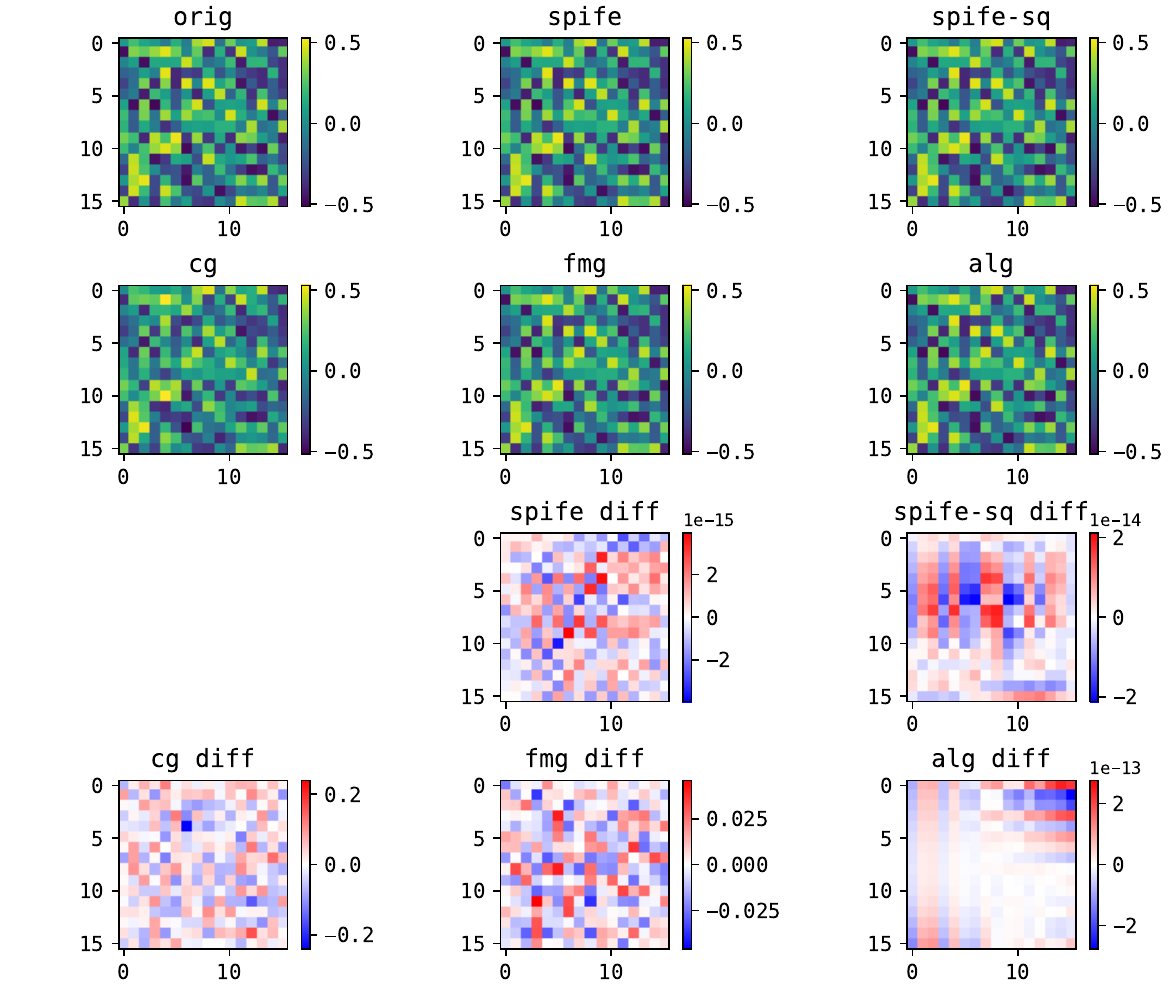}
    \end{minipage}
    \\
    \hline
    (b)
    &
    \begin{minipage}{0.6\textwidth}
    \includegraphics[width=1.0\textwidth]{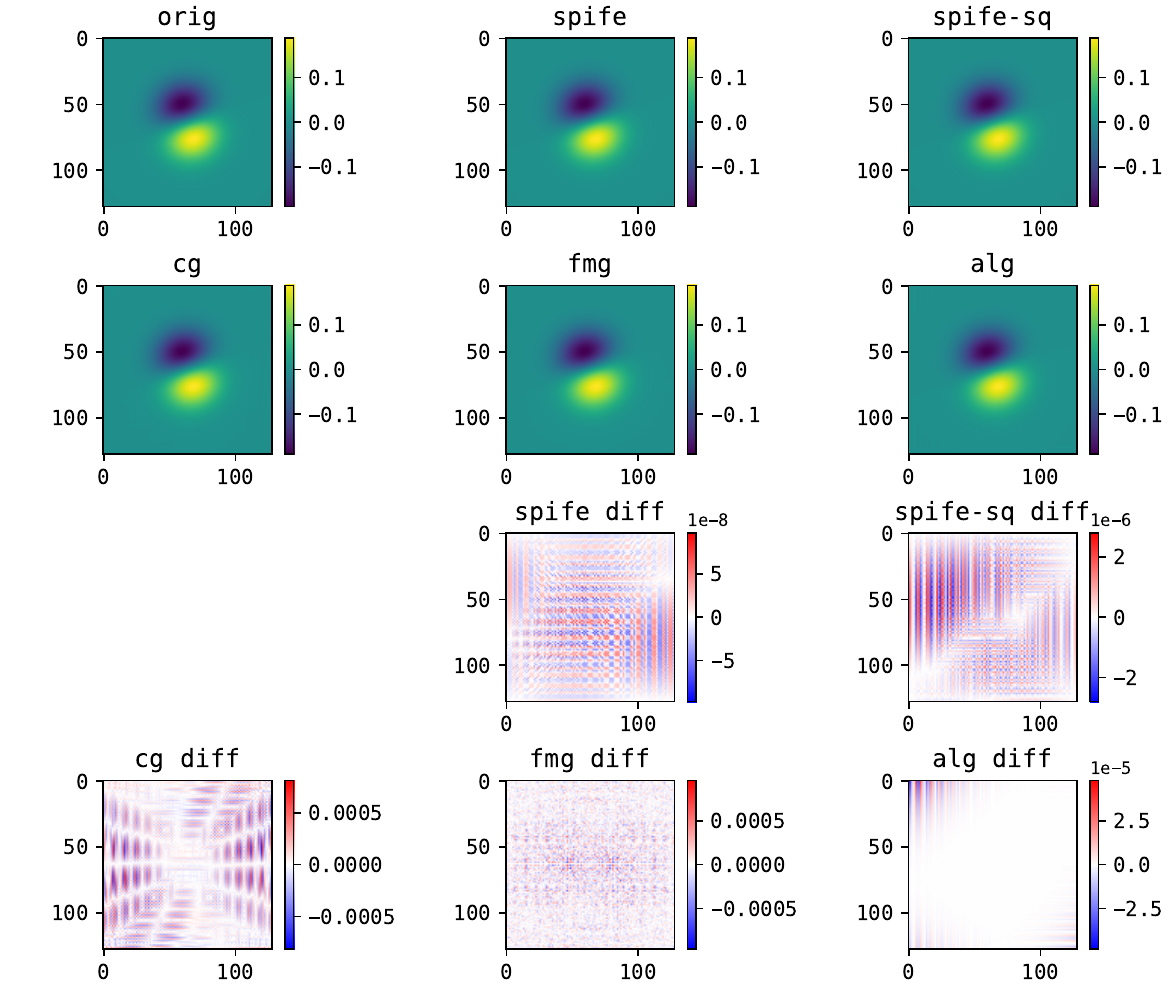}
    \end{minipage}
    \\
    \hline
    (c)
    &
    \begin{minipage}{0.6\textwidth}
    \includegraphics[width=1.0\textwidth]{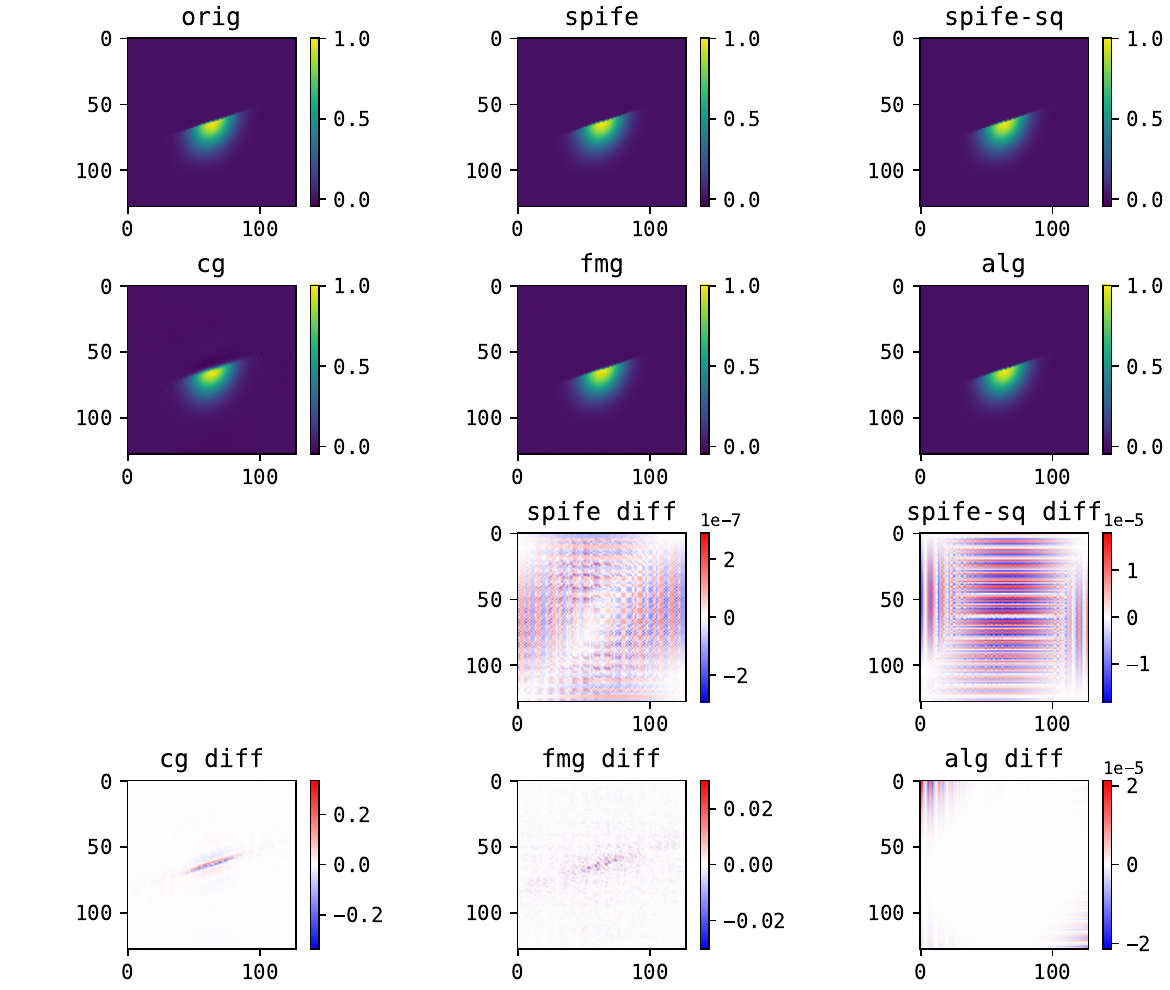}
    \end{minipage}
    \\
    \hline
  \end{tabular}
  \caption{Comparison of inversion results for (a) a random image
  of size $N=16$, 
  (b) smooth wave packet pulse of size $N=128$, and (c) truncated Gaussian of size $N=128$.
  For all iterative methods, $\log N$ iterations were used.
  }
  \label{fig:inv}
\end{figure}

We try inverting various images from their ADRT data using these methods. The
resulting inverse and the difference with the original image are shown in
Fig.~\ref{fig:inv}.
Three images are tested: a random $16 \times 16$ image where each pixel is
drawn from $\cU[-\half, \half]$, a smooth $128 \times 128$ image of
a sinusoidal wave packet with Gaussian amplitude, and 
a $128 \times 128$ image of a truncated Gaussian. 

In all of the tested images, the error for {\tt spife} is below {\tt cg}
and {\tt fmg}: {\tt 1e-15} for $16 \times 16$ images and below {\tt 1e-7} for
$128 \times 128$ images. {\tt spife-sq} and {\tt alg} behave similarly but have
higher errors of a couple more orders of magnitude.
{\tt fmg} and {\tt cg} have errors around {\tt 1e-1} and {\tt 1e-2} for $16
\times 16$ and {\tt 1e-4} and {\tt 1e-5} for $128 \times 128$. In both cases,
{\tt fmg} has errors roughly an order of magnitude lower in comparison.

Next, we add different types of noise to the ADRT data above and observe
how the accuracy of the inversion result changes. The noise perturbation causes
the pseudo-inverse to deviate from the image because now the ADRT data is no
longer in $\rg(\nR)$. Throughout these experiments, the accuracy of {\tt cg} and
{\tt fmg} inverses is not affected; this is to be expected since these
iterative schemes contain a projection to the range (these results also confirm
that the noise perturbations are small enough so that the solution to the normal
equations is still close to the original image).  

\begin{figure}
  \centering
  \begin{tabular}{r|c}
    \hline\\
    (a)
    &
    \begin{minipage}{0.6\textwidth}
    \includegraphics[width=1.0\textwidth]{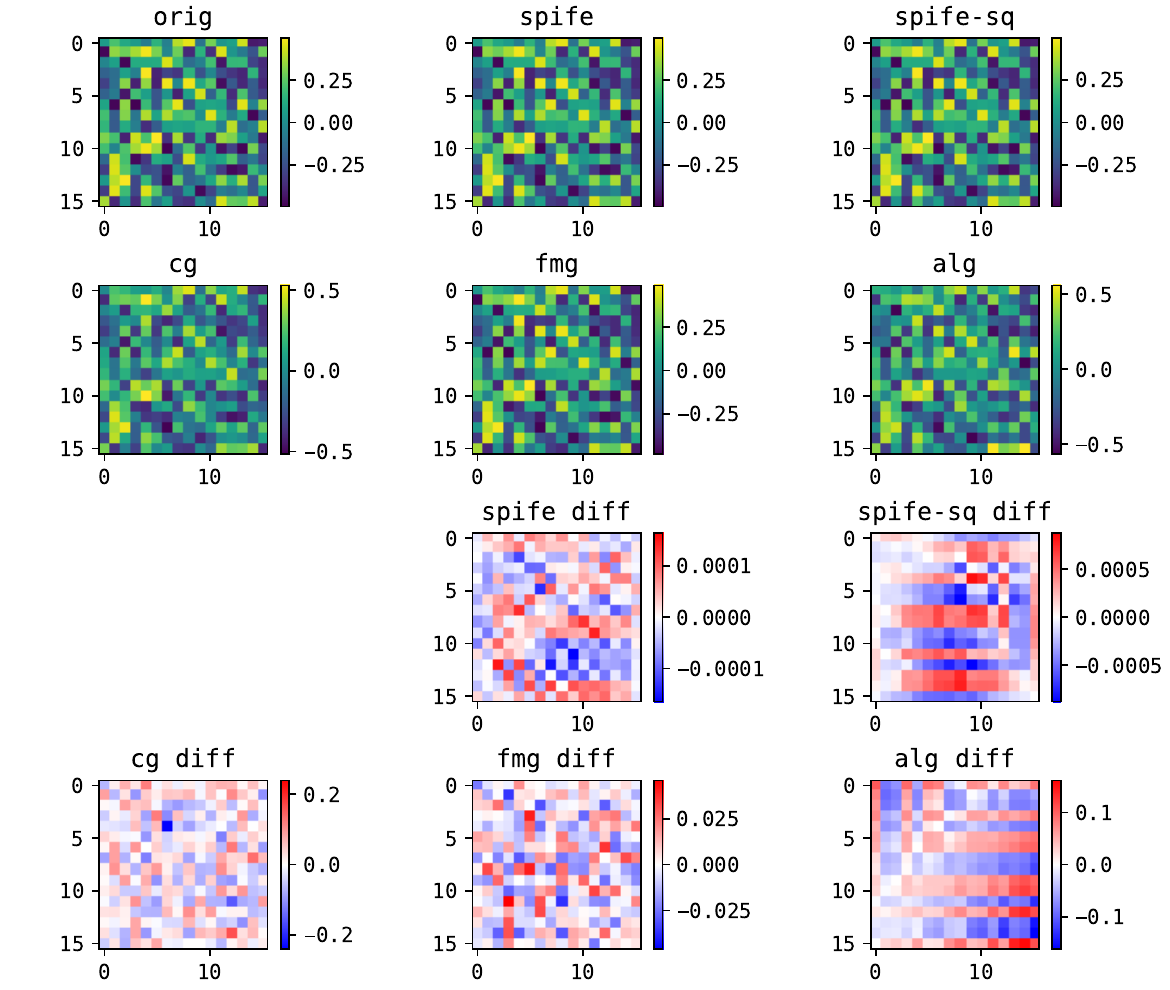}
    \end{minipage}
  \\ \hline
    (b)
    &
    \begin{minipage}{0.6\textwidth}
    \includegraphics[width=1.0\textwidth]{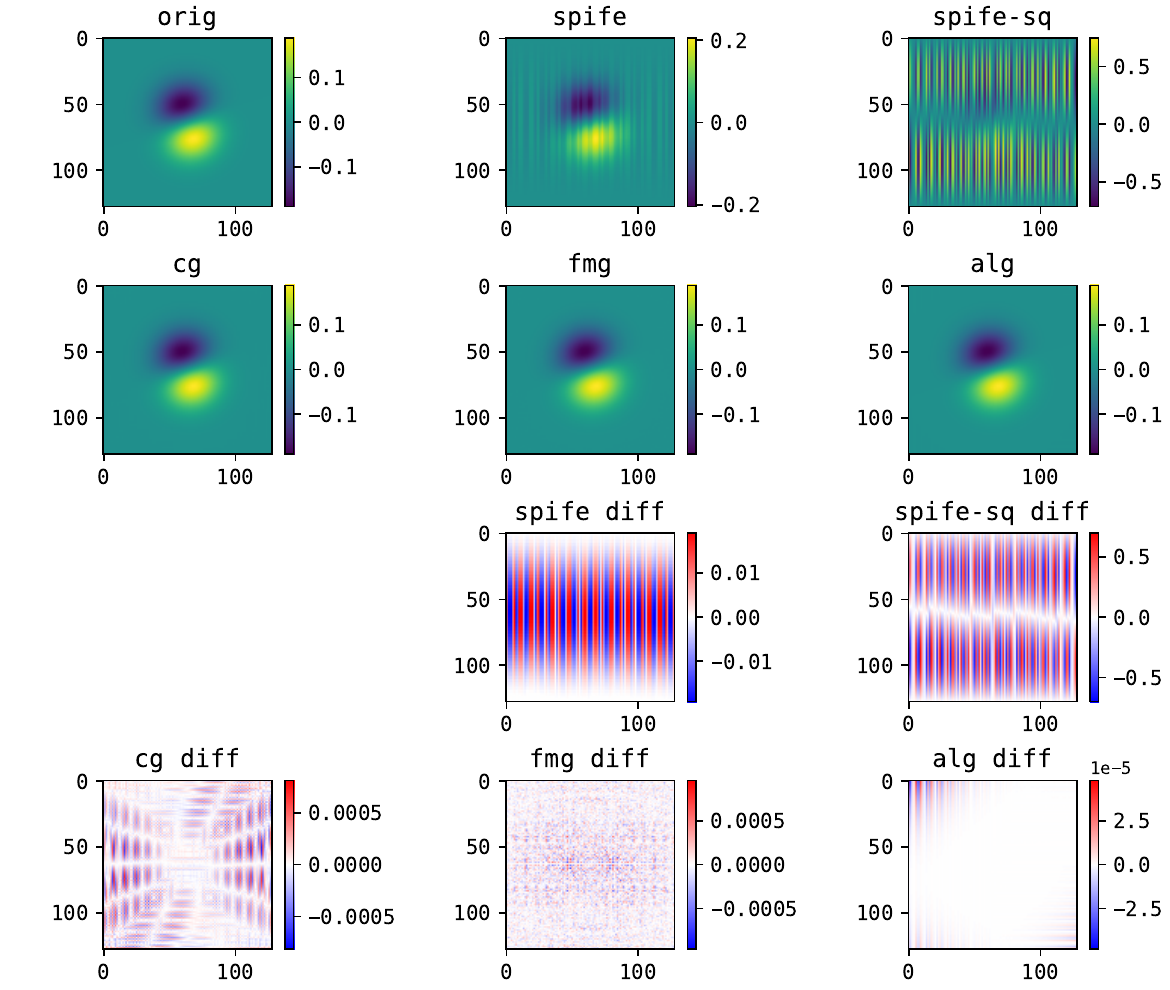}
    \end{minipage}
  \\ \hline
    (c)
    &
    \begin{minipage}{0.6\textwidth}
  \includegraphics[width=1.0\textwidth]{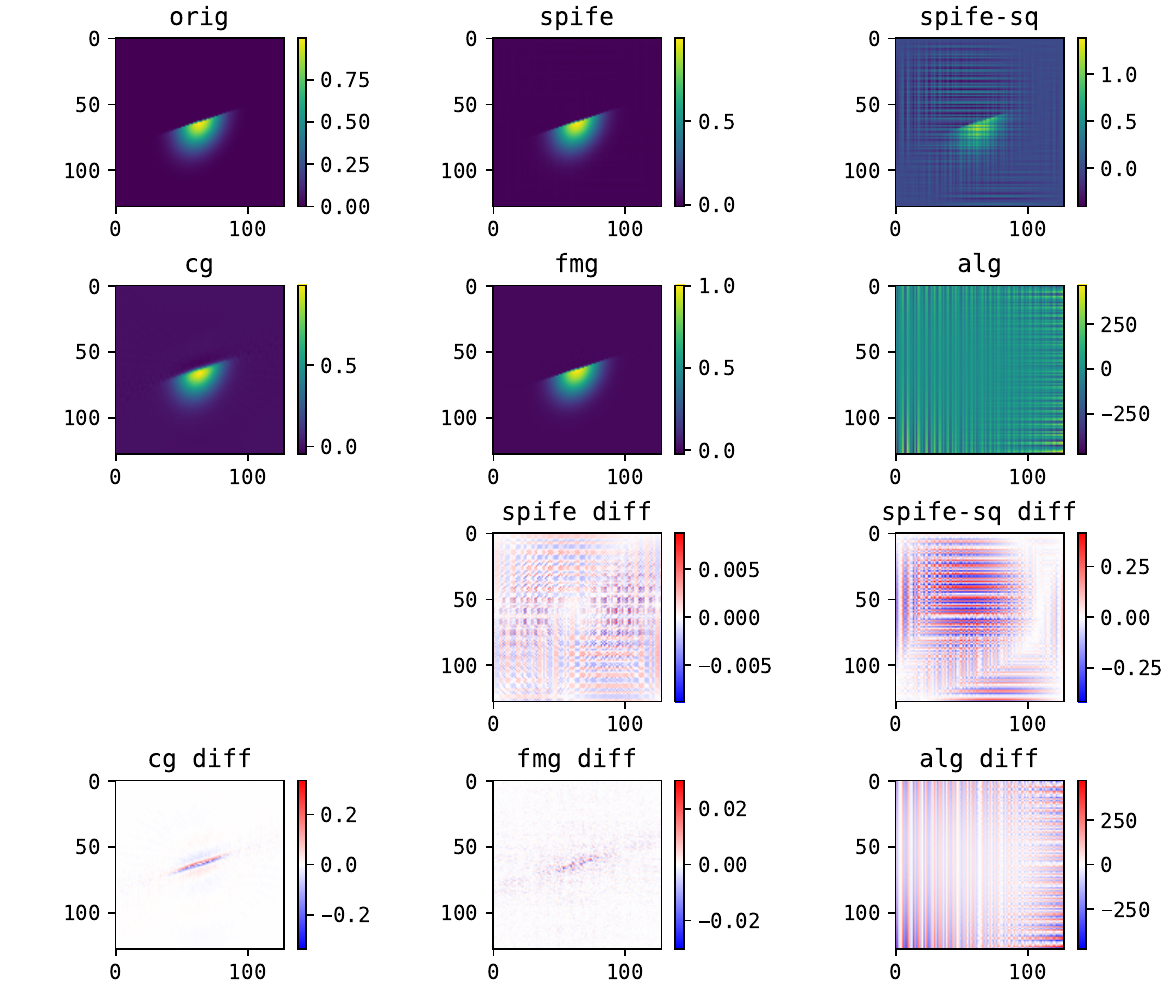}
    \end{minipage}
  \\
  \hline
  \end{tabular}
  \caption{Inversion results on the perturbed ADRT data. Results when 
  (a) uncorrelated uniform random noise with $\sigma=${\tt 2e-3} (SNR
  0.1\%) was added to the data with $N = 16$, (b) a single pixel was
  perturbed to amplitude {\tt 1e-7} in the data with $N=128$ (SNR {\tt
  2e-18}\%), and (c) uncorrelated random noise from a normal distribution
  with $\sigma=${\tt 1e-5} was added to the data with $N=128$ (SNR {\tt
  4e-20}\%). For all iterative methods, $\log N$ iterations were used.}
  \label{fig:add_noise}
\end{figure}

 The noise levels were chosen to result in similar orders of inversion
error for {\tt spife}. Since the stability of SPIFE depends on the image
size, different levels of noise were used for different $N$.
We add uncorrelated uniform random noise of amplitude $0.01$ for the random
image used in Fig.~\ref{fig:inv}(a). 
The results are shown in Fig.~\ref{fig:add_noise}(a).
{\tt alg} shows
dramatic loss of accuracy; the errors render the image unrecognizable. The {\tt
spife} inverse contains about {\tt 1e-4} error, performing better than {\tt
spife-sq} with errors around {\tt 5e-4}. This is attributable to the fact that
cross-quadrant level $m=1$ decomposition (Thm.~\ref{thm:spd}) tends to project
out the unstable modes, with a stabilizing effect.

Next, we perturb the smooth wave packet pulse image with $N=128$ used in Fig.~\ref{fig:inv} (b) with normal noise with $\sigma =${\tt1e-5}.
The results are shown in Fig.~\ref{fig:add_noise}(b).
In this case, the behavior is similar to that of the uniform noise.
The difference between {\tt spife} and {\tt spife-sq} is more dramatic,
due to the additional growth of an instability through $7$ single-quadrant
levels versus $4$. 

Finally, we perturb the truncated Gaussian image with $N=128$ from Fig.~\ref{fig:inv}(c) by
adding to a single pixel in the ADRT data by a value
of {\tt 1e-7}.
The results are shown in Fig.~\ref{fig:add_noise}(c).
This example highlights the difference between {\tt alg} and the
new variants {\tt spife} and {\tt spife-sq}. The perturbed pixel does not affect
the algebraically exact inverse, since the exact formula relies on local
calculations and the perturbed pixel in this case is positioned in such a way that
the other pixels are not affected globally. For the spectral pseudo-inverses,
however, the local perturbation effect is global, exciting high-frequency modes
for the quadrant where the perturbed pixel was placed.

\begin{figure}
  \centering
  \includegraphics[width=1.0\textwidth]{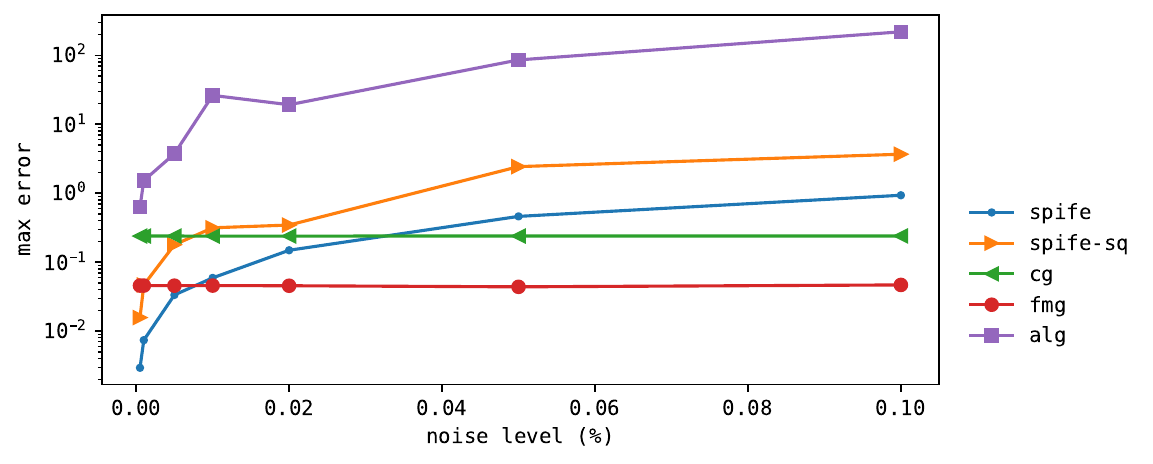}
  \caption{Maximum error in the inverse for a given noise level.}
  \label{fig:noise_levels}
\end{figure}

Taking the randomized image in the previous experiments (used in the top row of
Fig.~\ref{fig:add_noise}), 
we perform a test for various noise levels. The result is
shown in Fig.~\ref{fig:noise_levels}: {\tt spife} is overall more accurate compared
to the other non-iterative inverses.

\begin{figure}[h!]
  \includegraphics[width=0.9\textwidth]{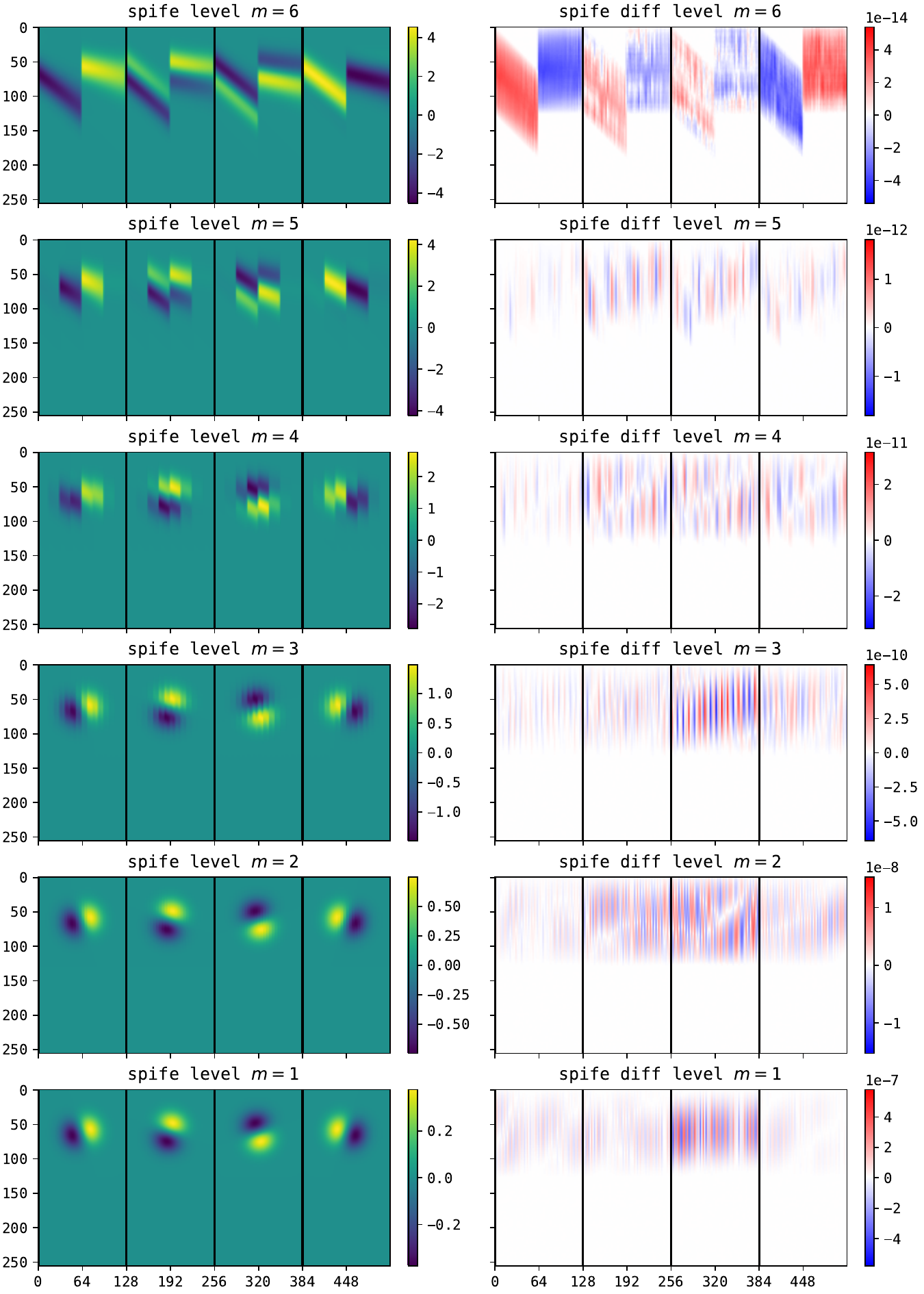}
  \caption{Growth of numerical error during the intermediate levels during the
  pseudo-inverse {\tt spife} (Thm.~\ref{prop:spifefull}).}
  \label{fig:error_growth}
\end{figure}

Finally, we illustrate the growth of numerical error during the repeated
pseudo-inverses throughout the levels for the smooth image of higher resolution
$256 \times 256$, as shown in Fig.~\ref{fig:error_growth}. 

\begin{figure}
  \includegraphics[width=0.8\textwidth]{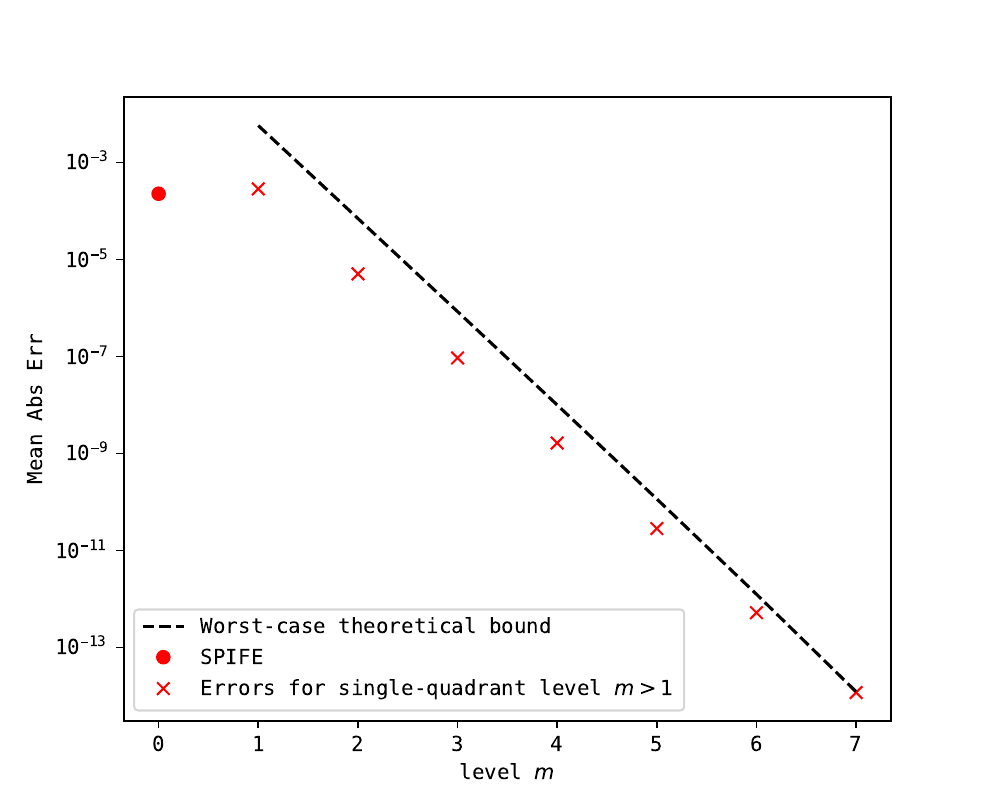}
  \caption{The growth of error during the single-quadrant levels
  for a random image with each pixel value sampled from the standard normal 
  distribution. The theoretical worst-case error bound is given in
  \eqref{eq:error_bound}. }
  \label{fig:error_growth_more}
\end{figure}

The rate of growth in error during the single-quadrant levels is upper
bounded by the product of the smallest singular values. Denoting by $\cE_m$ the 
$L^1$ error of the SPIFE inverse at level $m$, we may estimate based 
on \eqref{eq:Ktsk},
\begin{equation} \label{eq:error_bound}
  \sigma_{\text{min}, m}^+ \cE_{m} \le  \cE_{m+1},
  \where
  \sigma_{\text{min}, m}^+
  :=
  2 \cos \left( \frac{\pi}{2} - \frac{1}{2^n + 2^{m} + 1} \right).
\end{equation}
This straightforward bound is pessimistic and over-estimates the error seen in
practice.
In Fig.~\ref{fig:error_growth_more}, $\cE_m$ are shown alongside the theoretical
bound \eqref{eq:error_bound} for a randomly sampled image of size $128 \times
128$. The error is seen to grow at a steady rate during the single-quadrant
levels, but the error decreases moderately upon the projection to the
range of the cross-quadrant level $m=1$.

In broad terms, the spectral pseudo-inverse stabilizes the exact
inversion formula, although the effect of data perturbation becomes global. For
moderately small image sizes $N \le 2^8$, the method yields superior accuracy
when compared to both the CG iterations and the full multigrid method, run with
comparable computational cost. The growth of numerical errors in larger images
suggests a need for a stabilization technique for the pseudo-inverse.

\section{Conclusion} We have derived a new spectral decomposition of
the ADRT that yields an explicit formula called SPIFE that inverts the ADRT with
optimal complexity $\cO(N^2 \log^2 N)$ (Thm.~\ref{prop:spifefull}). We showed
through numerical examples that SPIFE is a competitive algorithm in terms of
accuracy vs. workload for in-range ADRT data of a moderately-sized image.
However, we also showed that SPIFE is highly sensitive to perturbations outside
the range of the ADRT. We illustrated this by deriving a simple estimate
for the error growth during the single-quadrant levels ($m > 1$) and by
numerical examples. 

The in-range assumption is typically not satisfied in applications, but the
spectral decomposition is still of interest as it suggests new research
directions: The decomposition can be used to
precondition existing iterative approaches; if the range characterization of
the ADRT can be further developed to yield a fast projection algorithm, SPIFE can
potentially be modified to yield an explicit inversion algorithm that is
accurate for general data not necessarily in range. These are interesting
open problems that will be considered in the near future.

\section*{Acknowledgements} 
WL is supported by NSF-DMS Award {\tt \#2309602}, a PSC-CUNY grant, and a
start-up fund from the Foundation for City College.
Part of this work was done while DR was enjoying
the hospitality of the International Research Institute of Disaster Science
(IRIDeS) at Tohoku University supported by Invitational Fellowship Program for
Collaborative Research with International Researcher (FY2023) hosted by Kenjiro
Terada.
The work of KR is partially supported by the NSF-DMS Award {\tt \#1937254}.

% spbasic: Springer medical, life sciences, chemistry, geology, engineering and computer science publications
% spmpsci.bst: Springer mathematics, computer science, and physical sciences journals publications.   
% spphys.bst: Springer physics publications
% \bibliographystyle{spmpsci}
\bibliographystyle{siamplain}

\end{document}